\documentclass[12pt]{article}
\usepackage{amsmath}
\usepackage{amssymb, amsbsy, amstext}
\usepackage{srcltx}
\usepackage{amscd,amsfonts,color}
\usepackage{float}
\usepackage{cite}

\usepackage{graphicx}
\usepackage{stmaryrd}
\usepackage{mathabx}
\usepackage{enumerate}
\usepackage{mathrsfs}
\usepackage{tikz}
\usepackage{booktabs}
\usepackage{caption}
\usepackage{pgflibraryarrows}
\usepackage{pgflibrarysnakes}

\input amssym.def
\input amssym.tex

\newtheorem{lem}{Lemma}[section]%
\newtheorem{theo}[lem]{Theorem}%
\newtheorem{defi}[lem]{Definition}%
\newtheorem{cor}[lem]{Corollary}%
\newtheorem{exam}[lem]{Example}%
\newtheorem{prob}[lem]{Open Problem}%
\newtheorem{prop}[lem]{Proposition}%
\newtheorem{rem}[lem]{Remark}%
\newtheorem{alm}[lem]{Algorithm}

\newcommand{\Aut}{\hbox{\rm Aut\,}}

\newcommand{\Ker}{\hbox{\rm Ker\,}}

\newcommand{\Core}{\hbox{\rm Core}}

\newcommand{\demo}{\f {\bf Proof}\hskip10pt}

\newcommand{\M}{\mathcal{M}}
\newcommand{\HH}{\mathcal{H}}

\newcommand{\med}{\hbox{\rm med}}
\newcommand{\dig}{\hbox{\rm dig}}

\def\Sym{\hbox{\rm Sym}}
\def\a{\alpha} \def\b{\beta}

    \def\G{\varGamma}
\def\ol1{\overline 1}
\def\di{\bigm|}
\def\lg{\langle}
\def\rg{\rangle}

\def\f{\noindent}
\newcommand{\qed}{\hfill \mbox{\raisebox{0.7ex}{\fbox{}}}}

\begin{document}

\title{
Linear hypermaps--modelling linear hypergraphs on surfaces\thanks{
{\it E-mail address}: pktide@163.com (Kai Yuan), wqgzyx23@163.com (Qi Wang), fengrq@math.pku.edu.cn (Rongquan Feng), wang$_{-}$yan@pku.org.cn(Yan Wang).}
}

\date{}

\author{Kai Yuan$^1$, Qi Wang$^1$,  Rongquan Feng$^{2}$, Yan Wang$^1$\\
{\small$^1$ School of Mathematics and Information Science, Yantai University, Yantai 264005, China} \\[-0.8ex]
{\small$^2$ School of Mathematical Sciences, Peking University, Beijing 100871, China }}

\maketitle

\begin{abstract}
A hypergraph is {\it linear} if each pair of distinct vertices appears in at most one common edge.  We say $\G=(V,E)$ is an {\it associated graph} of a linear hypergraph $\HH=(V, X)$ if for any $x\in X$, the induced subgraph $\G[x]$ is a cycle, and for any $e\in E$, there exists a unique edge $y\in X$ such that $e\subseteq y$. A {\it linear hypermap} $\mathcal{M}$ is a $2$-cell embedding of a connected linear hypergraph $\mathcal{H}$'s associated graph $\G$ on a compact connected surface, such that for any edge $x\in E(\HH)$, $\G[x]$ is the boundary of a $2$-cell and for any $e\in E(\G)$, $e$ is incident with two distinct $2$-cells. In this paper, we introduce linear hypermaps to model linear hypergraphs on surfaces and regular linear hypermaps modelling configurations on the surfaces. As an application, we classify regular linear hypermaps on the sphere and  determine the total number of proper regular linear hypermaps of genus 2 to 101.
\end{abstract}

{\small {\em 2000 Mathematics Subject Classification}:\,
              05C25, 05C30, 05C65, 05E15, 20B25.

             }

\section{Introduction}

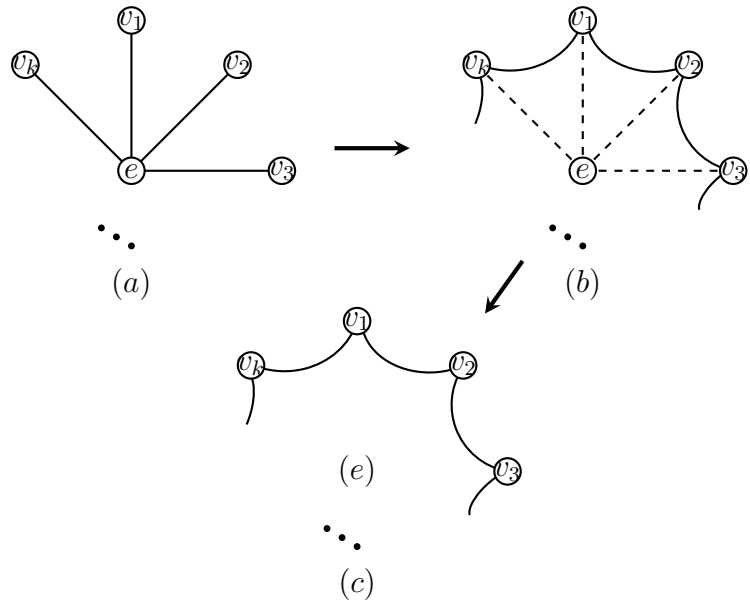
\begin{figure}[H]
\begin{center}
\begin{tikzpicture}

\draw[thick](0,0)--(2,0);
\draw[thick](0,0)--(0,2);
\draw[thick](0,0)--(1.41,1.41);
\draw[thick](0,0)--(-1.41,1.41);

\filldraw[fill=white, draw=black][thick](0,0)circle(5pt);
\filldraw[fill=white, draw=black][thick](2,0)circle(5pt);
\filldraw[fill=white, draw=black][thick](0,2)circle(5pt);
\filldraw[fill=white, draw=black][thick](1.41,1.41)circle(5pt);
\filldraw[fill=white, draw=black][thick](-1.41,1.41)circle(5pt);

\node at (0,0) {$e$};
\node at (2,0) {$v_3$};
\node at (0,2) {$v_1$};
\node at (1.41,1.41) {$v_2$};
\node at (-1.41,1.41) {$v_k$};

\fill (0,-1) circle (1.3pt);
\fill (-0.2,-0.88) circle (1.3pt);
\fill (-0.4,-0.76) circle (1.3pt);
\node at (0,-1.15) [below] {$(a)$};

\tikzstyle{arrow} = [ultra thick,->,>=stealth]
\draw[arrow] (2.7,0.3) -- (3.7,0.3);

\draw[dashed, thick](6,0)--(8,0);
\draw[dashed, thick](6,0)--(6,2);
\draw[dashed, thick](6,0)--(7.41,1.41);
\draw[dashed, thick](6,0)--(4.59,1.41);
\draw[thick] (7.41,1.41) arc (300:170:0.9 and 0.7);
\draw[thick] (8,0) arc (260:150:0.9 and 0.9);
\draw[thick] (6,2) arc (345:250:1 and 0.9);
\draw[thick,rotate=75] (2.45,-4.1) arc (350:250:0.5 and 0.2);
\draw[thick,rotate=220] (-5.45,5.25) arc (350:250:0.5 and 0.2);

\filldraw[fill=white, draw=black][thick](6,0)circle(5pt);
\filldraw[fill=white, draw=black][thick](8,0)circle(5pt);
\filldraw[fill=white, draw=black][thick](6,2)circle(5pt);
\filldraw[fill=white, draw=black][thick](7.41,1.41)circle(5pt);
\filldraw[fill=white, draw=black][thick](4.59,1.41)circle(5pt);

\node at (6,0) {$e$};
\node at (8,0) {$v_3$};
\node at (6,2) {$v_1$};
\node at (7.41,1.41) {$v_2$};
\node at (4.59,1.41) {$v_k$};

\fill (6,-1) circle (1.3pt);
\fill (5.8,-0.88) circle (1.3pt);
\fill (5.6,-0.76) circle (1.3pt);
\node at (6,-1.15) [below] {$(b)$};

\draw[arrow] (5.2,-1.2) -- (4.7,-1.9);

\draw[thick] (4.41,-2.59) arc (300:170:0.9 and 0.7);
\draw[thick] (5,-4) arc (260:150:0.9 and 0.9);
\draw[thick] (3,-2) arc (345:250:1 and 0.9);
\draw[thick,rotate=75] (-2.2,-2.2) arc (350:250:0.5 and 0.2);
\draw[thick,rotate=220] (-0.5,6.4) arc (350:250:0.5 and 0.2);

\filldraw[fill=white, draw=black][thick](5,-4)circle(5pt);
\filldraw[fill=white, draw=black][thick](3,-2)circle(5pt);
\filldraw[fill=white, draw=black][thick](4.41,-2.59)circle(5pt);
\filldraw[fill=white, draw=black][thick](1.59,-2.59)circle(5pt);

\node at (3,-4) {$(e)$};
\node at (5,-4) {$v_3$};
\node at (3,-2) {$v_1$};
\node at (4.41,-2.59) {$v_2$};
\node at (1.59,-2.59) {$v_k$};

\fill (3,-5) circle (1.3pt);
\fill (2.8,-4.88) circle (1.3pt);
\fill (2.6,-4.76) circle (1.3pt);
\node at (3,-5.15) [below] {$(c)$};

\end{tikzpicture}
\caption{The modification process.}
\label{modif}
\end{center}
\end{figure}

Let $\HH$ be a hypergraph  with vertex set $V(\HH)$ and hyperedge set $ E(\HH)$. The associated {\it Levi graph} $\mathcal{G}$ of $\HH$ is a bipartite graph with
vertex set $V(\HH) \sqcup E(\HH)$, in which $v \in V(\HH)$ and $e \in E(\HH)$
are adjacent if and only if $v\in e$ in $\HH$. Given a 2-cell embedding of  $\mathcal{G}$ into a compact and connected surface $\mathcal{S}$ without border, we modify this embedding to obtain an embedding of
$\HH$  into $\mathcal{S}$ wherein  certain regions  of the modified embedding  represent  hyperedges of $\HH$, and the remaining regions are the   faces of this embedding.
 As in \cite{Wh1984,Wh}, we illustrate the modification process in Figure \ref{modif}.  In part (a) of
Figure \ref{modif}, we see - in the embedding of $\mathcal{G}$ - the vertex $e$ which is an hyperedge consisting of vertices $ \{v_1,v_2,\ldots,v_k\}$ of $\HH$. In part (b), we modify the imbedding of $\mathcal{G}$ by adding edges $\{v_i,v_{i+1}\}\    (\mod k)$,
 deleting the vertex $e$ and the edges $\{e, v_i\}, 1 \le i \le k$, so that
in part (c) the hyperedge $e$ stands for a region in the modified
embedding.

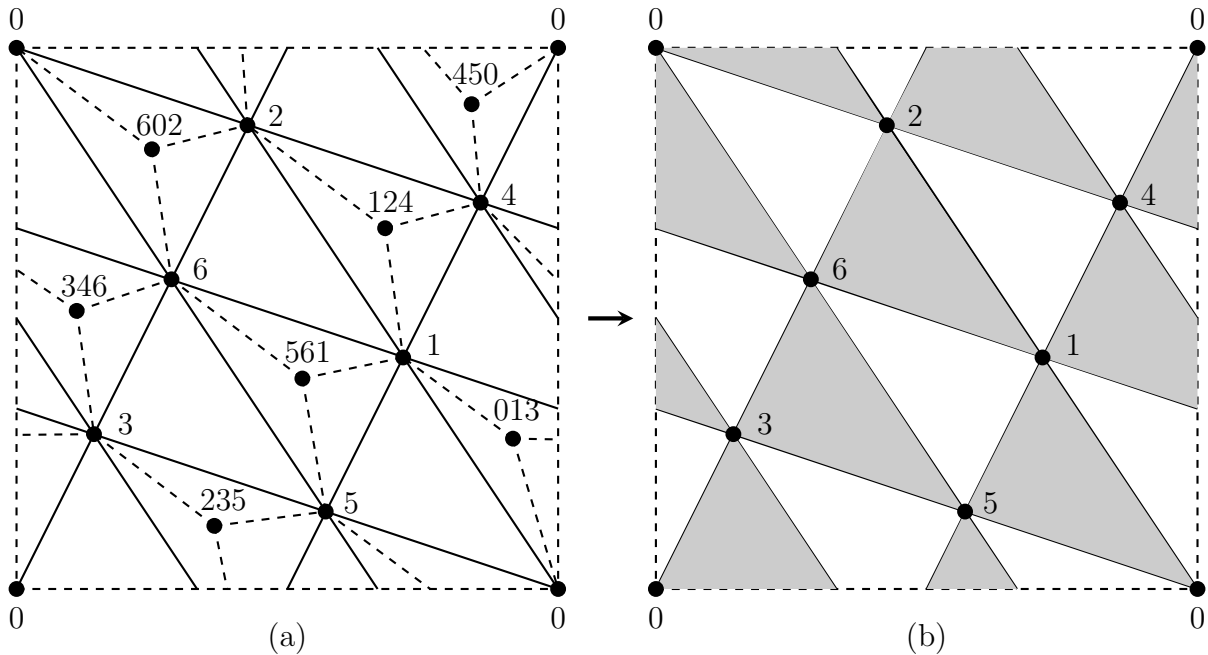
\begin{figure}[H]
\begin{center}
\begin{tikzpicture}
\begin{scope}[xshift=-4cm]

\node at (0,0.1) [above] {0};
\node at (7.2,0.1) [above] {0};
\node at (0,-7.3) [below] {0};
\node at (7.2,-7.3) [below] {0};

\draw [dashed, thick] (0,0)--(7.2,0);
\draw [dashed, thick] (0,0)--(0,-7.2);
\draw [dashed, thick] (0,-7.2)--(7.2,-7.2);
\draw [dashed, thick] (7.2,0)--(7.2,-7.2);

\draw [thick] (0,0)--(7.2,-2.4);
\draw [thick] (0,-2.4)--(7.2,-4.8);
\draw [thick] (0,-4.8)--(7.2,-7.2);

\draw [thick] (3.6,0)--(0,-7.2);
\draw [thick] (7.2,0)--(3.6,-7.2);

\draw [thick] (4.8,0)--(7.2,-3.6);
\draw [thick] (2.4,0)--(7.2,-7.2);
\draw [thick] (0,0)--(4.8,-7.2);
\draw [thick] (0,-3.6)--(2.4,-7.2);

\fill (5.14,-4.12) circle (3pt);
\node at (5.3,-3.95) [right] {1};
\fill (3.07,-1.03) circle (3pt);
\node at (3.2,-0.9) [right] {2};
\fill (1.03,-5.14) circle (3pt);
\node at (1.2,-5) [right] {3};
\fill (6.17,-2.06) circle (3pt);
\node at (6.3,-1.95) [right] {4};
\fill (4.11,-6.17) circle (3pt);
\node at (4.2,-6.05) [right] {5};
\fill (2.06,-3.08) circle (3pt);
\node at (2.2,-2.94) [right] {6};

\draw [dashed, thick] (3.07,-1.03)--(3,0);
\draw [dashed, thick] (1.03,-5.14)--(0,-5.15);
\draw [dashed, thick] (6.17,-2.06)--(7.2,-3.1);
\draw [dashed, thick] (4.11,-6.17)--(5.5,-7.2);

\fill (1.8,-1.35) circle (3pt);
\node at (1.9,-1.3) [above] {602};
\draw [dashed, thick] (1.8,-1.35)--(0,0);
\draw [dashed, thick] (1.8,-1.35)--(3.07,-1.03);
\draw [dashed, thick] (1.8,-1.35)--(2.06,-3.08);

\fill (0.8,-3.5) circle (3pt);
\node at (0.9,-3.45) [above] {346};
\draw [dashed, thick] (0.8,-3.5)--(0,-2.95);
\draw [dashed, thick] (0.8,-3.5)--(1.03,-5.14);
\draw [dashed, thick] (0.8,-3.5)--(2.06,-3.08);

\fill (6.05,-0.75) circle (3pt);
\node at (6.1,-0.63) [above] {450};
\draw [dashed, thick] (6.05,-0.75)--(7.2,0);
\draw [dashed, thick] (6.05,-0.75)--(6.17,-2.06);
\draw [dashed, thick] (6.05,-0.75)--(5.4,0);

\fill (4.9,-2.4) circle (3pt);
\node at (4.97,-2.33) [above] {124};
\draw [dashed, thick] (4.9,-2.4)--(3.07,-1.03);
\draw [dashed, thick] (4.9,-2.4)--(6.17,-2.06);
\draw [dashed, thick] (4.9,-2.4)--(5.14,-4.12);

\fill (3.8,-4.4) circle (3pt);
\node at (3.87,-4.33) [above] {561};
\draw [dashed, thick] (3.8,-4.4)--(2.06,-3.08);
\draw [dashed, thick] (3.8,-4.4)--(5.14,-4.12);
\draw [dashed, thick] (3.8,-4.4)--(4.11,-6.17);

\fill (2.63,-6.36) circle (3pt);
\node at (2.75,-6.3) [above] {235};
\draw [dashed, thick] (2.63,-6.36)--(1.03,-5.14);
\draw [dashed, thick] (2.63,-6.36)--(4.11,-6.17);
\draw [dashed, thick] (2.63,-6.36)--(2.8,-7.2);

\fill (6.6,-5.2) circle (3pt);
\node at (6.65,-5.1) [above] {013};
\draw [dashed, thick] (6.6,-5.2)--(5.14,-4.12);
\draw [dashed, thick] (6.6,-5.2)--(7.2,-5.21);
\draw [dashed, thick] (6.6,-5.2)--(7.2,-7.2);

\fill (0,0) circle (3pt);
\fill (7.2,0) circle (3pt);
\fill (0,-7.2) circle (3pt);
\fill (7.2,-7.2) circle (3pt);

\node at (3.6,-7.5) [below] {(a)};

\tikzstyle{arrow} = [ultra thick,->,>=stealth]
\draw[arrow] (7.6,-3.6) -- (8.2,-3.6);
\end{scope}
\begin{scope}[xshift=4.5cm]

\draw [dashed, thick] (0,0)--(7.2,0);
\draw [dashed, thick] (0,0)--(0,-7.2);
\draw [dashed, thick] (0,-7.2)--(7.2,-7.2);
\draw [dashed, thick] (7.2,0)--(7.2,-7.2);

\draw [thick] (0,0)--(7.2,-2.4);
\draw [thick] (0,-2.4)--(7.2,-4.8);
\draw [thick] (0,-4.8)--(7.2,-7.2);

\draw [thick] (3.6,0)--(0,-7.2);
\draw [thick] (7.2,0)--(3.6,-7.2);

\draw [thick] (4.8,0)--(7.2,-3.6);
\draw [thick] (2.4,0)--(7.2,-7.2);
\draw [thick] (0,0)--(4.8,-7.2);
\draw [thick] (0,-3.6)--(2.4,-7.2);

\fill[gray!40] (3.07,-1.03)--(0,0)--(2.4,0)--(3.07,-1.03);
\fill[gray!40] (2.06,-3.08)--(0,0)--(0,-2.4)--(2.06,-3.08);
\fill[gray!40] (1.03,-5.14)--(0,-4.8)--(0,-3.6)--(1.03,-5.14);

\fill[gray!40] (3.07,-1.03)--(3.6,0)--(4.8,0)--(6.17,-2.06)--(3.07,-1.03);
\fill[gray!40] (2.06,-3.08)--(5.14,-4.12)--(3.07,-1.03)--(2.06,-3.08);
\fill[gray!40] (2.06,-3.08)--(1.03,-5.14)--(4.11,-6.17)--(2.06,-3.08);
\fill[gray!40] (1.03,-5.14)--(0,-7.2)--(2.4,-7.2)--(1.03,-5.14);

\fill[gray!40] (6.17,-2.06)--(7.2,0)--(7.2,-2.4)--(6.17,-2.06);
\fill[gray!40] (6.17,-2.06)--(5.14,-4.12)--(7.2,-4.8)--(7.2,-3.6)--(6.17,-2.06);
\fill[gray!40] (4.11,-6.17)--(5.14,-4.12)--(7.2,-7.2)--(4.11,-6.17);
\fill[gray!40] (4.11,-6.17)--(4.8,-7.2)--(3.6,-7.2)--(4.11,-6.17);

\fill (5.14,-4.12) circle (3pt);
\node at (5.3,-3.95) [right] {1};
\fill (3.07,-1.03) circle (3pt);
\node at (3.2,-0.9) [right] {2};
\fill (1.03,-5.14) circle (3pt);
\node at (1.2,-5) [right] {3};
\fill (6.17,-2.06) circle (3pt);
\node at (6.3,-1.95) [right] {4};
\fill (4.11,-6.17) circle (3pt);
\node at (4.2,-6.05) [right] {5};
\fill (2.06,-3.08) circle (3pt);
\node at (2.2,-2.94) [right] {6};

\fill (0,0) circle (3pt);
\fill (7.2,0) circle (3pt);
\fill (0,-7.2) circle (3pt);
\fill (7.2,-7.2) circle (3pt);

\node at (0,0.1) [above] {0};
\node at (7.2,0.1) [above] {0};
\node at (0,-7.3) [below] {0};
\node at (7.2,-7.3) [below] {0};

\node at (3.6,-7.5) [below] {(b)};
\end{scope}
\end{tikzpicture}

\caption{The Fano plane on Torus.}
\label{Fano}

\end{center}
\end{figure}

In Figure \ref{Fano}, we illustrate the modification process of embedding the Fano plane $\HH$ on Torus, where $E(\HH) = \{\{0, 1, 3\}, \{1, 2, 4\}, $ $\{2, 3, 5\},
 \{3, 4, 6\}, \{4, 5, 0\}, \{5, 6, 1\}, \{6, 0, 2\}\}$.
Part (a) of Figure \ref{Fano} shows the embedding of the Levi graph $\mathcal{G}$ together with the added edges (the solid lines) around each hyperedge. Part (b) is the modified embedding of $\HH$, where the shaded regions are the faces and the white ``regions" depict the hyperedges of $\HH$, see \cite{Wh1995,FC,Wh2002}.

A hypergraph is {\it linear} if each pair of distinct vertices appear in at most one common hyperedge.  A {\em configuration} $(v_r,b_k)$ is a finite incidence structure with $v$ points and $b$ blocks, with the property that (i) each block contains exactly $k$ points; (ii) each point is contained in exactly $r$ blocks; and (iii) any pair of distinct points is contained in at most one block. Clearly, a configuration is a linear hypergraph. Hypergraph embedding has been described by embedding of the associated Levi graph, called the hypermaps, see \cite{BF,BF1,BJ, C, Co,CPS,CM,DH,DY,WY,J,JM,Ra}.

Finite geometries can be modeled by graph imbeddings, see \cite{EJS,Wh1995,FC}.  To model  finite geometries on surfaces, Arthur T. White said in \cite{Wh2002}   that the challenge is to find an imbedding of particular interest: maximum characteristic, maximal symmetry, revealing properties of the geometry, and
so forth.  This paper is motivated by modelling configurations on surfaces with maximal symmetry \cite{Wh2002} and we primarily achieve our expectations by employing methods rooted in group theory.

\begin{defi}
A graph $\G=(V,E)$ is called an {\it associated graph} of a linear hypergraph $\HH=(V, X)$ if for any $x\in X$, the induced subgraph $\G[x]$ is a cycle, and for any $e\in E$ with two end vertices $u_e, v_e$, there exists a unique hyperedge $y\in X$ such that $\{u_e, v_e\}\subseteq y$.
\end{defi}
Specially, if~$|x|=1$, then~$\G[x]$ is a loop; and if~$|x|=2$, then~$\G[x]$ is a connected graph with two multiple edges. For each $x\in X$, let $\widehat{\G}[x]$ be a cycle. Let $\widehat{\G}$ be the union of graphs $\widehat{\G}[x]$ for $x\in X$. Clearly, $\widehat{\G}$ is an associated graph of $\HH$ and the associated graph of $\HH$ may not be unique. A linear hypergraph is said to be {\it loopless} if its associated graph is without loops.
A linear hypermap can be defined in a combinatorial way as follows.

\begin{defi}
Suppose that $\HH$ is a connected linear hypergraph. A linear hypermap $\M$ is a $2$-cell embedding of a graph $\G$ which is an associated graph of $\HH$  on a compact connected surface $\mathcal{S}$, such that for any hyperedge $x\in E(\HH)$, $\G[x]$ is the boundary of a $2$-cell and for any $e\in E(\G)$, $e$ is incident with exactly two distinct $2$-cells.
\end{defi}

The genus of $\M$ is defined as the genus of $\mathcal{S}$. A $2$-cell of $\M$ is called a {\it hyperedge} if it is bounded by $\G[x]$ for some $x\in E(\HH)$; otherwise, it is called a {\it hyperface}. The hypergraph $\HH$ is called the {\it underlying linear hypergraph} of $\M$ and the graph $\G$ is called the {\it associated graph} of $\M$.

 In this paper, in order to investigate linear hypermaps with high symmetry, we always make the following assumption.

 \vspace{0.3cm}
 {\noindent \bf Assumption.}

  \vspace{0.2cm}

 {\noindent $\bullet$ The underlying linear hypergraph of a linear hypermap has more than one hyperedge and loopless.}

  \vspace{0.3cm}

\begin{figure}[H]
\begin{center}
\begin{tikzpicture}

\draw[dashed, thick](-2.3,-2.7)--(-2.3,4.1)--(6.4,4.1)--(6.4,-2.7)--(-2.3,-2.7);
\fill[gray!30](-2.3,-2.7)--(-2.3,4.1)--(6.4,4.1)--(6.4,-2.7)--(-2.3,-2.7);

\filldraw[fill=white, draw=black][thick](0,0)circle(1.6);
\filldraw[fill=white, draw=black][thick](3.6,0)circle(2);
\filldraw[fill=white, draw=black][thick](1.5,1.9)circle(0.8);
\filldraw[fill=white, draw=black][thick](3.3,2.75)circle(0.75);
\filldraw[fill=white, draw=black][thick](2.2,2.8)circle(0.35);
\filldraw[fill=white, draw=black][thick](-0.8,1.8)circle(0.35);
\filldraw[fill=white, draw=black][thick](1.58,-1)circle(0.25);

\end{tikzpicture}

\caption{7 coins on the Sphere.}
\label{coin}
\end{center}
\end{figure}
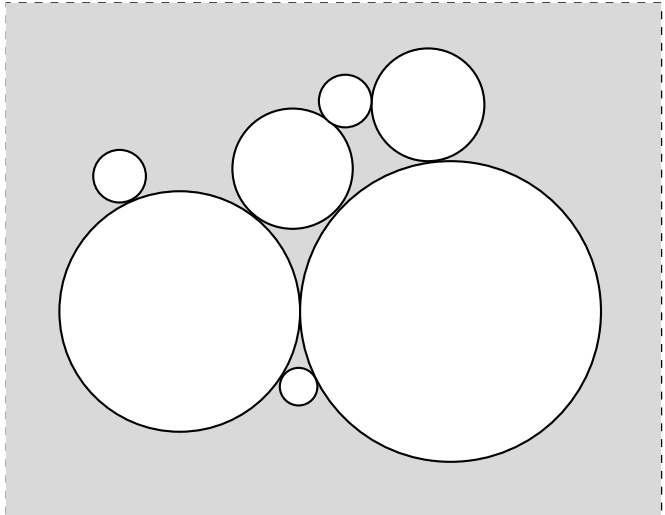

If $\mathcal{S}$ is orientable, then $\M$ is called orientable; otherwise, non-orientable. Clearly, one can look the linear hypermap $\M$ as a map $M$ on $\mathcal{S}$ whose vertices remain the same and whose faces are the hyperedges and hyperfaces. In this way, $M$ is called the {\it associated map} of $\M$ which is clearly the embedding of the associated graph $\G$. So a linear hypermap is a map representation of a linear hypergraph. Linear hypermaps also have connections with two-orbit polyhedra \cite{IH}, and with digraphs embedding on orientable surfaces whose anti-faces are hyperedges and faces are hyperfaces \cite{BCM}.   Moreover, we can use a linear hypermap to model coins (white components) on surfaces, see Figure \ref{coin}. In this case, every vertex (the point of tangency between two circles) is incident with 2 hyperedges (white components). Therefore, linear hypermaps are intriguing.

\begin{figure}[H]
\begin{center}
\begin{tikzpicture}
\pgfmathsetmacro{\n}{6}
\pgfmathsetmacro{\radius}{1}
\foreach \i in {1,...,\n} {
\coordinate (p\i) at ({360/\n * (\i - 1)}:\radius);}
\fill[gray!50] (p1)--(p2)--(p3)--(p4)--(p5)--(p6)--(p1);
\node at (p1) [left] {4};
\node at (p2) [below] {3};
\node at (p3) [below] {2};
\node at (p4) [right] {9};
\node at (p5) [above] {8};
\node at (p6) [above] {6};
\foreach \i in {6,...,2} {\draw[thick] (p\i) -- (p\the\numexpr\i-1\relax);}
\draw[thick] (p1) -- (p6);

\pgfmathsetmacro{\n}{6}
\pgfmathsetmacro{\radius}{2}
\foreach \i in {1,...,\n} {
\coordinate (t\i) at ({360/\n * (\i - 1)}:\radius);}
\foreach \i in {6,...,2} {\draw[dashed] (t\i) -- (t\the\numexpr\i-1\relax);}
\draw[dashed] (t1) -- (t6);

\pgfmathsetmacro{\n}{12}
\pgfmathsetmacro{\radius}{1.72}
\foreach \i in {1,...,\n} {
\coordinate (s\i) at ({360/\n * (\i - 1)}:\radius);}
\foreach \i in {1,...,6} {\fill (s\the\numexpr\i*2\relax) circle (1pt);}
\foreach \i in {2} {\node at (s\i) [right] {5};}
\foreach \i in {4} {\node at (s\i) [above] {1};}
\foreach \i in {6} {\node at (s\i) [left] {7};}
\foreach \i in {8} {\node at (s\i) [left] {5};}
\foreach \i in {10} {\node at (s\i) [below] {1};}
\foreach \i in {12} {\node at (s\i) [right] {7};}
\draw[thick] (s2)--(p2);
\draw[thick] (p1)--(s2);
\draw[thick] (s4)--(p3);
\draw[thick] (p2)--(s4);
\draw[thick] (s6)--(p4);
\draw[thick] (p3)--(s6);
\draw[thick] (s8)--(p5);
\draw[thick] (p4)--(s8);
\draw[thick] (s10)--(p6);
\draw[thick] (p5)--(s10);
\draw[thick] (s12)--(p1);
\draw[thick] (p6)--(s12);

\fill[gray!50] (s2)--(p2)--(s4)--(t2)--(s2);
\fill[gray!50] (s4)--(p3)--(s6)--(t3)--(s4);
\fill[gray!50] (s6)--(p4)--(s8)--(t4)--(s6);
\fill[gray!50] (s8)--(p5)--(s10)--(t5)--(s8);
\fill[gray!50] (s10)--(p6)--(s12)--(t6)--(s10);
\fill[gray!50] (s12)--(p1)--(s2)--(t1)--(s12);

\foreach \i in {1,...,6} {\fill (p\i) circle (1.4pt);}
\foreach \i in {1,...,6} {\fill (s\the\numexpr\i*2\relax) circle (1.4pt);}
\end{tikzpicture}

\caption{A linear hypermap on Torus.}
\label{torus}
\end{center}
\end{figure}
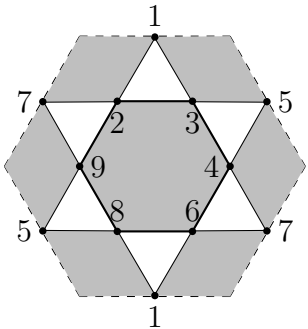

Our example in Figure \ref{torus} is a linear hypermap on  Torus with 9 vertices, 6 hyperedges (white components) and 3 hyperfaces (grey components). The hyperedges of the underlying linear hypergraph are $X=\{\{1,2,3\},\{3,4,5\},\{4,6,7\},\{1,6,8\},\{5,8,9\},\{2,7,9\}\}$.

\section{Algebraic Linear Hypermaps}

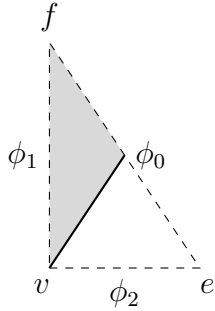
\begin{figure}[H]
\begin{center}
\begin{tikzpicture}

\fill[gray!30](0,0)--(1,1.5)--(0,3)--(0,0);
\draw[dashed](0,0)--(2,0)--(0,3)--(0,0);
\draw[thick](0,0)--(1,1.5);
\node at (-0.1,0) [below] {$v$};
\node at (2.1,0) [below] {$e$};
\node at (0,3) [above] {$f$};
\node at (0,1.5) [left] {$\phi_1$};
\node at (1,1.5) [right] {$\phi_0$};
\node at (1,0) [below] {$\phi_2$};

\end{tikzpicture}
\caption{\it A flag of a linear hypermap.}
\label{flag}
\end{center}
\end{figure}

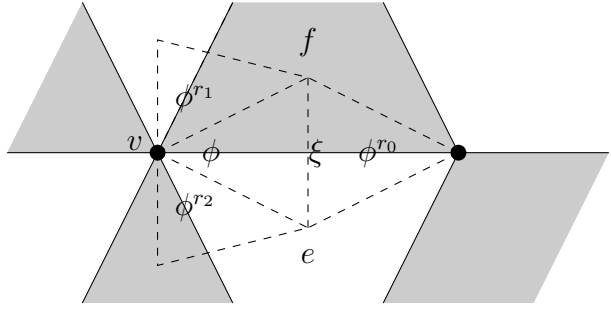
\begin{figure}[H]
\begin{center}
\begin{tikzpicture}[
dashed line/.style={dashed, thin},
dot/.style = {
      draw,
      fill = white,
      circle,
      inner sep = 0pt,
      minimum size = 5pt
    }
]
\fill[gray!40] (0,0)--(4,0)--(3,2)--(1,2)--(0,0);
\fill[gray!40] (0,0)--(1,-2)--(-1,-2)--(0,0);
\fill[gray!40] (0,0)--(-2,0)--(-1,2)--(0,0);
\fill[gray!40] (4,0)--(6,0)--(5,-2)--(3,-2)--(4,0);
\draw (4,0) -- (3,2)   (0,0)--(4,0) (0,0) -- (1,2) (0,0) -- (-1,2) (0,0) -- (-1,-2) (0,0) -- (-2,0) (0,0) -- (1,-2) (4,0) -- (3,-2) (4,0)--(6,0);
\draw[dashed](0,0) -- (1,2) (0,0) -- (0,1.5)--(2,1) (0,0) -- (0,-1.5)--(2,-1) (0,0) -- (2,1) (0,0) -- (2,-1) (2,1) -- (2,-1)   (4,0) -- (2,1) (4,0) -- (2,-1) ;
\fill (0,0) circle (3pt);
\fill (4,0) circle (3pt);
\path (0.3,0) node[label = {right: $\phi$}] {};
\path (2.5,0) node[label = {left: $\xi$}] {};
\path (3.5,0) node[label = {left: $\phi^{r_0}$}] {};
\path (0.5,1.2) node[label = {below: $\phi^{r_1}$}] {};
\path (0.5,-1.2) node[label = {above: $\phi^{r_2}$}] {};
\path (2,1) node[label = {above: $f$}] {};
\path (2,-1) node[label = {below: $e$}] {};
\path (-0.3,0.5) node[label = {below: $v$}] {};
\end{tikzpicture}

\caption{\it $r_i$ acting on $\phi = (\phi_0, \phi_1, \phi_2)$ for $i\in \{0,1,2\}$.}
\label{ahm}
\end{center}
\end{figure}

Linear hypermaps can also be defined in an algebraic way.
A flag $\phi$, see Figure \ref{flag} or Figure \ref{ahm}, of a linear hypermap is a triangle with two colored parts (divided by an edge of its associated graph), three vertices labeled $v,e$ and $f$, and three sides labeled $\phi_0, \phi_1$ and $\phi_2$. The white and the grey parts of $\phi$ belong to a hyperedge and a hyperface, respectively; the three vertices $v, e, f$ represent a triad of mutually incident vertex, hyperedge, and hyperface; the three sides $\phi_0, \phi_1$ and $\phi_2$ are opposite to  $v,e$ and $f$, respectively. In this paper, the flag $\phi$ is denoted by $\phi=(\phi_0, \phi_1, \phi_2; \xi)$, where $\phi_0, \phi_1, \phi_2$ are called the $0$-dimensional,  $1$-dimensional and $2$-dimensional components of $\phi$, respectively, and $\xi$ is the unique edge contained in it.
\begin{rem}

{\rm (1)}
It is obvious that a flag of a linear hypermap is one to one corresponding to an arc of its associated graph. So, for brevity in this paper, we usually denote a flag by $\phi=(\phi_0, \phi_1, \phi_2)$ without listing the arc when we emphasize on the hypermap but not on the graph.

{\rm (2)} A flag of a linear hypermap contains two flags of its associated map which is divided by the inner edge.
\end{rem}

Let $\M$ be a linear hypermap and let $\Phi$ be the set of flags of $\M$. Take $\phi=(\phi_0, \phi_1, \phi_2) \in \Phi$, then for each $i\in \{0, 1, 2\}$, there exists a unique flag $\hat{\phi}=(\hat{\phi}_0, \hat{\phi}_1, \hat{\phi}_2)$ sharing the same $i$-dimensional item with $\phi$. Define three involutions on $\Phi$ which preserve the incidence relations of flags, denoted by $r_i, i\in \{0,1,2\}$, such that
\begin{enumerate}[{\rm(i)}]
\item If $\hat{\phi}_0=\phi_0$, then $\phi^{r_0}=\hat{\phi}$ and $\hat{\phi}^{r_0}=\phi; $
\item If $\hat{\phi}_1=\phi_1$, then $\phi^{r_1}=\hat{\phi}$ and $\hat{\phi}^{r_1}=\phi;$
\item If $\hat{\phi}_2=\phi_2$, then $\phi^{r_2}=\hat{\phi}$ and $\hat{\phi}^{r_2}=\phi$.
\end{enumerate}
 Let $G= \langle r_0, r_1, r_2\rangle$, then $G$ is called the monodromy group of $\M$. Obviously, $G$ is a subgroup of  $S_\Phi$ and $G$ acts transitively on $\Phi$ if the underlying hypergraph of $\M$ is connected.
The following Figure \ref{ahm} shows how $r_0, r_1$ and $ r_2$ act on a typical flag $\phi=(\phi_0, \phi_1, \phi_2)$.

\begin{rem}
The three involutions $r_0,r_1$ and $r_2$  are pairwise distinct and fixed-point-free for the underlying linear hypergraph is connected with more than one hyperedge and loopless.
\end{rem}
\begin{prop}\label{inc}
Let $\phi,\hat{\phi}\in \Phi$. Then,
\begin{enumerate}[{\rm(1)}]
\item $\phi$ and $\hat{\phi}$ are incident with the common vertex if and only if $\hat{\phi}\in \phi^{\lg r_1,r_2 \rg};$
\item $\phi$ and $\hat{\phi}$ are incident with the common hyperedge if and only if $\hat{\phi}\in \phi^{\lg r_0,r_2 \rg};$
\item $\phi$ and $\hat{\phi}$ are incident with the common hyperface if and only if $\hat{\phi}\in \phi^{\lg r_0,r_1 \rg}$.
\end{enumerate}
\end{prop}

\begin{prop}
$\langle r_1,r_2\rangle\cap \langle r_0,r_2\rangle=\langle r_2\rangle$.
\end{prop}
\demo
Let $\phi=(\phi_0, \phi_1, \phi_2)\in \Phi$ and $H=\langle r_1,r_2\rangle$. Clearly, $\phi^{\langle r_1,r_2\rangle}\cap \phi^{\langle r_0,r_2\rangle}=\phi^{\langle r_2\rangle}.$ Let $a\in \langle r_1,r_2\rangle\cap \langle r_0,r_2\rangle$. It suffices to show that $a\in \langle r_2\rangle$. Obviously,
$
\phi^a\in \phi^{\langle r_2\rangle}
 $ and it is suitable for any $\phi\in \Phi$.
 It follows that $\phi^a$ and $\phi$ are incident with the common hyperedge and $(\phi^{r_2})^a\in \phi^{r_2\langle r_2\rangle}=\phi^{\langle r_2\rangle}$.

 Because $\phi^{\langle r_2\rangle}=\{\phi, \phi^{r_2}\}$, there are two possibilities: $\phi^a=\phi, (\phi^{r_2})^a=\phi^{r_2}$ or $\phi^a=\phi^{r_2}, (\phi^{r_2})^a=\phi$. In either case, one can get $ar_2=r_2a$ easily.
So, $a\in C_H(\langle r_2\rangle)$. If $|r_1r_2|$ is odd, then $a\in C_H(\langle r_2\rangle)=\langle r_2\rangle$, as required. If $|r_1r_2|=n$ is even, then $C_H(\langle r_2\rangle)=\langle (r_1r_2)^{\frac{n}{2}}, r_2\rangle \cong \mathbb{Z}_2\times \mathbb{Z}_2$. Suppose that $a\in \langle r_2\rangle(r_1r_2)^{\frac{n}{2}}$, then there exits $\tilde{\phi}$ such that $\tilde{\phi}$ and $\tilde{\phi}^a$ are incident with different hyperedges, contradicting to $\tilde{\phi}^a \in \tilde{\phi}^{\langle r_2\rangle}$ which implies that $\tilde{\phi}^a$ and $\tilde{\phi}$ are incident with the same hyperedge. Therefore, $a\in \langle r_2\rangle$.
\qed

\begin{prop}
Let $\phi=(\phi_0, \phi_1, \phi_2)\in \Phi$. Then $\phi^{\langle r_1,r_2\rangle\langle r_0,r_2\rangle}\cap \phi^{\langle r_0,r_2\rangle\langle r_1,r_2\rangle}=\phi^{\langle r_1,r_2\rangle\cup\langle r_0,r_2\rangle}.$
\end{prop}
\demo Obviously,
$$\phi^{\langle r_1,r_2\rangle\cup\langle r_0,r_2\rangle} \subseteq  \phi^{\langle r_1,r_2\rangle\langle r_0,r_2\rangle}\cap \phi^{\langle r_0,r_2\rangle\langle r_1,r_2\rangle}.$$
Let $\phi$ be incident with the vertex $v$, the hyperedge $e$, and the hyperface $f$. Assume that there are $a,d\in\langle r_1,r_2\rangle$ and $b,c\in\langle r_0,r_2\rangle$ such that $\phi^{ab}=\phi^{cd}$.
It suffices to show that $\phi^{ab}\in \phi^{\langle r_1,r_2\rangle\cup\langle r_0,r_2\rangle}$. Set $\hat{\phi}=\phi^{ab}=(\hat{\phi}_0, \hat{\phi}_1, \hat{\phi}_2)=\phi^{cd}$ which is incident with the vertex $v^\prime$, the hyperedge $e^\prime$, and the hyperface $f^\prime$. By Proposition \ref{inc}, one can obtain that $v=v^\prime$, or $e=e^\prime$, or $v\ne v^\prime$ and $e\ne e^\prime$.
If $v=v^\prime$, then $\phi^{ab}\in \phi^{\langle r_1,r_2\rangle}$;
if $e=e^\prime$, then $\phi^{ab}\in \phi^{\langle r_0,r_2\rangle}$.
Now suppose that $v\ne v^\prime$ and $e\ne e^\prime$. Then $\phi^a$ is incident with $v$ and $e^\prime$, and $\phi^c$ is incident with $v^\prime$ and $e$.
It follows that both $v$ and $v^\prime$ are incident with $e$ and $e^\prime$, which contradicts to the assumption of $\M$ being a linear hypermap. Therefore $\phi^{\langle r_1,r_2\rangle\langle r_0,r_2\rangle}\cap \phi^{\langle r_0,r_2\rangle\langle r_1,r_2\rangle}=\phi^{\langle r_1,r_2\rangle\cup\langle r_0,r_2\rangle},$ as desired.
\qed

\vskip 0.3cm

\begin{defi}\label{def2.5}
For a given finite set $\Phi$ and three pairwise distinct and fixed-point-free involutory permutations $r_0, r_1,r_2$ on $\Phi$, a quadruple $\mathcal{M}=\mathcal{M}(\Phi;r_0, r_1,r_2)$ is called an algebraic linear hypermap if the following hold:
\begin{enumerate}[{\rm(1)}]
\item $\langle r_1,r_2\rangle\cap \langle r_0,r_2\rangle=\langle r_2\rangle;$
\item $\phi^{\langle r_1,r_2\rangle\langle r_0,r_2\rangle}\cap \phi^{\langle r_0,r_2\rangle\langle r_1,r_2\rangle}=\phi^{\langle r_1,r_2\rangle\cup\langle r_0,r_2\rangle},$ for all $\phi\in \Phi;$
\item the monodromy group $\langle r_0, r_1, r_2\rangle$ acts transitively on $\Phi$.
\end{enumerate}
The orbits of  the subgroup $\langle r_1,r_2\rangle$ $($resp. $\langle r_0,r_2\rangle$
and $\langle r_0,r_1\rangle$$)$ are flags incident with vertices $($resp. hyperedges and hyperfaces$)$.
\end{defi}

\begin{rem}
Let $\phi \in \Phi$. Condition $(1)$ of Definition~\ref{def2.5} ensures that no hyperedges contain duplicated vertices, and condition $(2)$ ensures the underlying hypergraph to be linear.
The orbits of $\langle r_1,r_2\rangle$ $($resp. $\langle r_0,r_2\rangle$
and $\langle r_0,r_1\rangle$$)$ are in one-to-one correspondence with the vertices $($resp. hyperedges and hyperfaces$)$.
Vertices, hyperedges and hyperfaces are mutually incident by non-empty intersection of the corresponding orbits.
\end{rem}

Given a combinatorial linear hypermap $\M$  with three fixed-point-free involutions acting on flags, we can construct an algebraic linear hypermap $\mathcal{M}(\Phi;r_0, r_1,r_2)$
 through the procedure described above. The transitivity follows easily from the connectivity of  the underlying linear hypergraph.
Conversely, given any  transitive permutation group $G=\langle r_2, r_1, r_0\rangle$ on a set $\Phi$, where $r_0, r_1,r_2$ are fixed-point-free involutions, one
can construct a combinatorial linear hypermap $\M$ by taking the vertices, hyperedges and hyperfaces as the orbits of the subgroups $\langle r_1,r_2\rangle$, $\langle r_0,r_2\rangle$
and $\langle r_0,r_1\rangle$ on $\Phi$, respectively, and the incidence relation respects  non-empty intersection.

\begin{theo}\label{3.2}
$(${\rm Euler formula for linear hypermaps}$)$ Let $\M=\M(\Phi;r_0, r_1,r_2)$
be a linear hypermap on a closed surface $\mathcal{S}$ of genus $g$ having $v$ vertices, $e$ hyperedges and $f$ hyperfaces. Then
\begin{enumerate}[{\rm(1)}]
\item $v + e + f - |\Phi|/2 = 2 - 2g$, if $\mathcal{S}$ is orientable$;$
\item $v + e + f - |\Phi|/2 = 2 - g$, if $\mathcal{S}$ is non-orientable.
\end{enumerate}
\end{theo}

\demo
Let $\G$ be the associated graph of $\M$. Then $\G$ has $v$ vertices, $|\Phi|/2$ edges and $e + f$ faces. The equalities follow from the Euler characteristic formula $\chi(\mathcal{S})=v  - |\Phi|/2  + e + f$.
\qed

\section{Linear Hypermaps Morphisms}

\begin{defi}
Let $\mathcal{M}=\mathcal{M}(\Phi;r_0, r_1,r_2)$ and $\widetilde{\mathcal{M}}=\mathcal{M}(\widetilde{\Phi};\tilde{r}_0, \tilde{r}_1,\tilde{r}_2)$ be two linear hypermaps. A morphism $\sigma$ from $\mathcal{M}$ to $\widetilde{\mathcal{M}}$ is a map $\sigma: \Phi\rightarrow \widetilde{\Phi}$ such that $ r_0\sigma=\sigma \tilde{r}_0, r_1\sigma=\sigma \tilde{r}_1$ and $r_2\sigma=\sigma \tilde{r}_2$.
\end{defi}

\begin{prop}\label{4.2}
Let $\mathcal{M}=\mathcal{M}(\Phi;r_0, r_1,r_2)$ and $\widetilde{\mathcal{M}}=\mathcal{M}(\widetilde{\Phi}; \tilde{r}_0, \tilde{r}_1, \tilde{r}_2)$ be two linear hypermaps with $G=\langle r_0, r_1, r_2 \rangle$ and $\widetilde{G}= \langle \tilde{r}_0, \tilde{r}_1, \tilde{r}_2 \rangle$. Let $\sigma: \mathcal{M}\rightarrow \widetilde{\mathcal{M}}$ be a linear hypermap morphism. Then $\sigma: \Phi\rightarrow \widetilde{\Phi}$ is surjective.
\end{prop}

\demo
Let $\tilde{\b}\in\widetilde{\Phi}$. We are required to show that $\tilde{\b}\in\Phi^\sigma$. Take some ~$\b\in\Phi^\sigma$, there exists $\a\in\Phi$ such that $\a^\sigma=\b$. Since $ \widetilde{G}$ acts transitively on $\widetilde{\Phi}$, there exists $\tilde{g}\in \widetilde{G}$ such that $\b^{\tilde{g}}=\tilde{\b}$. Assume that $\tilde{g}=\tilde{z_1}\tilde{z_2}\cdots \tilde{z_n}$ and $\tilde{z_i}\in \{\tilde{r}_0,\tilde{r}_1,\tilde{r}_2\},\, i\in \{1,2,\ldots,n\}$.
Set $g=z_1z_2\cdots z_n$, where $z_i\in \{r_0, r_1, r_2\}$ and $z_i=r_j$ whenever $\tilde{z_i}=\tilde{r_j}$.
 Since $r_0\sigma=\sigma \tilde{r}_0, r_1\sigma=\sigma \tilde{r}_1$ and $r_2\sigma=\sigma \tilde{r}_2$,
we have $\tilde{\b}=\b^{\tilde{g}}=\a^{\sigma \tilde{g}}=\a^{g\sigma}=(\a^g)^\sigma\in\Phi^\sigma$, as required.  \qed

\begin{prop}\label{4.3}
Let $\mathcal{M}=\mathcal{M}(\Phi;r_0, r_1,r_2)$ and $\widetilde{\mathcal{M}}=\mathcal{M}(\widetilde{\Phi}; \tilde{r}_0, \tilde{r}_1, \tilde{r}_2)$ be two linear hypermaps with $G=\langle r_0, r_1, r_2 \rangle$ and $\widetilde{G}= \langle \tilde{r}_0, \tilde{r}_1, \tilde{r}_2 \rangle$, and let $\sigma: \mathcal{M}\rightarrow \widetilde{\mathcal{M}}$ be a linear hypermap morphism. Then the following hold:
\begin{enumerate}[{\rm(1)}]
\item For any $g\in G$, there exists a unique $\tilde{g}\in \widetilde{G}$ such that $g\sigma=\sigma\tilde{g};$
\item Let $\tau$ be a map from $G$ to $\widetilde{G}$ defined as $g^\tau=\tilde{g}$ such that $g\sigma=\sigma\tilde{g}$ for any $g\in G$. Then $\tau$ is an epimorphism. In particular, if $\sigma$ is a bijection, then $\tau$ is an isomorphism.
\end{enumerate}
\end{prop}

\demo
(1) Since $r_0\sigma=\sigma \tilde{r}_0, r_1\sigma=\sigma \tilde{r}_1$, $r_2\sigma=\sigma \tilde{r}_2$, there exists $\tilde{g}\in \widetilde{G}$ such that $ g\sigma=\sigma \tilde{g}$. To show that $\tilde{g}$ is unique, it suffices to prove that if there exists $\tilde{g}_1\in \widetilde{G}$ such that $g\sigma=\sigma \tilde{g}_1$, then $\tilde{g}_1= \tilde{g}$. According to Proposition \ref{4.2},  for any $\b\in\widetilde{\Phi}$, there exists $\a\in\Phi$ such that $\b=\a^\sigma$. Hence, $\b^{\tilde{g}_1}=\a^{\sigma \tilde{g}_1}=\a^{\sigma \tilde{g}}=\b^{\tilde{g}}$. Therefore, $\tilde{g}_1= \tilde{g}$.

(2) For any two elements $x, y\in G$, by (1), $\sigma(xy)^\tau=(xy)\sigma=x\sigma y^\tau=\sigma(x^\tau y^\tau)$. So, $\tau$ is a group homomorphism. Since $r_0^\tau= \tilde{r}_0, r_1^\tau= \tilde{r}_1$ and $r_2^\tau= \tilde{r}_2$, $\tau$ is surjective. As a result, $\tau$ is an epimorphism.

Now assume that $\sigma$ is bijective. For any $\phi\in\Phi$, set $\phi^\sigma=\tilde{\phi}$, where $\tilde{\phi}\in\widetilde{\Phi}$.
Since $g\sigma=\sigma {g}^\tau$ and $\sigma$ is a bijection, it follows that $g\in G_\phi$ if and only if $  g^\tau\in\widetilde{G}_{\tilde{\phi}}$, so that ${G_\phi}^\tau=\widetilde{G}_{\tilde{\phi}}$. Since $\tau$ is a homomorphism, $(G_\phi\Ker \tau)^\tau={G_\phi}^\tau(\Ker \tau)^\tau=\widetilde{G}_{\tilde{\phi}}$. Thus $G_\phi\Ker \tau=G_\phi$,  that is, $\Ker \tau\subseteq G_\phi$. The arbitrary choice of
$\phi$ implies
$\Ker \tau\subseteq\bigcap_{\phi\in\Phi}G_\phi=1$, so $\Ker \tau=1$. Therefore, $\tau$ is injective and further implying $\tau$ to be an isomorphism.
\qed

\begin{defi}
Let $\mathcal{M}=\mathcal{M}(\Phi;r_0, r_1,r_2)$ and $\widetilde{\mathcal{M}}=\mathcal{M}(\widetilde{\Phi};\tilde{r}_0, \tilde{r}_1,\tilde{r}_2)$ be two linear hypermaps, and let $\sigma: \mathcal{M}\rightarrow \widetilde{\mathcal{M}}$ be a linear hypermap morphism.  If $\sigma$ is a bijection, then $\mathcal{M}$ is said to be isomorphic to $\widetilde{\mathcal{M}}$, denoted by $\mathcal{M}\cong \widetilde{\mathcal{M}}$. Especially, if $\mathcal{M}=\widetilde{\mathcal{M}}$, then $\sigma$ is said to be an automorphism of $\mathcal{M}$.
\end{defi}

 Let  $A = \Aut(\mathcal{M})$ be the {\it automorphism group} of $\mathcal{M}$, then $A=\{\sigma\di \sigma\in\Sym(\Phi), r_0\sigma=\sigma r_0, r_1\sigma=\sigma r_1, r_2\sigma=\sigma r_2\}$. The following Theorem~\ref{th3.5} about the automorphism group of linear hypermaps is trivial and we omit the proof.

\begin{theo}\label{th3.5}
 Let $\mathcal{M}=\mathcal{M}(\Phi;r_0, r_1,r_2), G=\langle r_0, r_1, r_2 \rangle$, and let $\Aut(\mathcal{M})=\{\sigma\di \sigma\in\Sym(\Phi), r_0\sigma=\sigma r_0, r_1\sigma=\sigma r_1, r_2\sigma=\sigma r_2\}$. Then,
$$\Aut(\mathcal{M})=C_{\Sym(\Phi)}(G)\cong N_{G}(H)/H,$$
where $H$ is the stabilizer of some $\phi \in \Phi$. Moreover, $\Aut(\mathcal{M})$ acts semi-regularly on $\Phi$.
\end{theo}

An automorphism of a linear hypermap $\M=\M(\Phi;r_0, r_1,r_2)$ can be regarded as a permutation of its vertices,  hyperedges and  hyperfaces that preserves incidence relations. Every automorphism of $\M$ is uniquely determined by its action on any incident vertex-hyperedge-hyperface triple $(v, e, f)$. So, the number of automorphisms of $\M$ is bounded above by the number of triples (which are sometimes called `flags'). A linear hypermap $\M$ is called regular if this upper bound is attained. That is, a linear hypermap ~$\M$ is {\it regular} if $\Aut(\M)$ acts transitively on $\Phi$.

\begin{theo}\label{4.8}
Let $\mathcal{M}=\mathcal{M}(\Phi;r_0, r_1,r_2)$ and $\widetilde{\mathcal{M}}=\mathcal{M}(\widetilde{\Phi};\tilde{r}_0, \tilde{r}_1,\tilde{r}_2)$ be two  regular  linear hypermaps with $G=\langle r_0, r_1,r_2 \rangle = \langle \tilde{r}_0, \tilde{r}_1, \tilde{r}_2 \rangle$. Then $\mathcal{M}\cong \mathcal{\widetilde{M}}$ if and only if there exists an automorphism $\tau$ of $G$ such that $r_0^\tau=\tilde{r}_0, r_1^\tau=\tilde{r}_1$ and $r_2^\tau=\tilde{r}_2$.
\end{theo}

\demo
 If $\M\cong \widetilde{\M}$, then there exists a bijection $\sigma: \Phi\rightarrow \widetilde{\Phi}$ such that $ r_0\sigma=\sigma \tilde{r}_0, r_1\sigma=\sigma \tilde{r}_1, r_2\sigma=\sigma \tilde{r}_2$.
By Proposition \ref{4.3}(2), there exists $\tau\in \Aut(G)$ such that $r_0^\tau=\tilde{r}_0, r_1^\tau=\tilde{r}_1, r_2^\tau=\tilde{r}_2$, as required.

 Conversely, take $\phi\in\Phi$ and $\tilde{\phi}\in\widetilde{\Phi}$. Set $\sigma : \phi^g\mapsto \widetilde{\phi}^{g^\tau},$ for all $g\in G$. Since $G$ and $\widetilde{G}$ act regularly on $\Phi$ and $\widetilde{\Phi}$, respectively,  it follows that $\sigma$ is bijective. For any $g\in G$ and $i\in \{0, 1, 2\}$, the following formulae hold: $$(\phi^g)^{r_i\sigma}=(\phi^{gr_i})^\sigma=\tilde{\phi}^{(gr_i)^\tau}
 =\tilde{\phi}^{g^\tau \tilde{r}_i}=
(\phi^g)^{\sigma\tilde{r}_i}.$$
As a result, $r_0\sigma=\sigma \tilde{r}_0, r_1\sigma=\sigma \tilde{r}_1$ and $r_2\sigma=\sigma \tilde{r}_2$. Hence, $\M\cong \widetilde{\M}$.
\qed

\vskip 0.3cm

Let $\mathcal{M}=\mathcal{M}(\Phi;r_0, r_1,r_2)$ be a regular linear hypermap. Since the centraliser of a transitive group must act semi-regularly, it follows that $G$ acts regularly and consequently
$$\Aut(\mathcal{M}) \cong G.$$
Therefore, the actions of $G$ and $\Aut(\mathcal{M})$ on $\Phi$ are permutations isomorphic to the right and the left regular representations of $G$, respectively. In this situation, we may identify $\Phi$ with $G$, and so $\mathcal{M}=\mathcal{M}(G;r_0, r_1,r_2)$. Therefore, a regular linear hypermap can be defined in the following way:

\begin{defi}
Let $G=\langle r_0, r_1, r_2\rangle$ be a finite group with three pairwise distinct involutory permutations $r_0, r_1,r_2$. A quadruple $\mathcal{M}=\mathcal{M}(G;r_0, r_1,r_2)$ is called a regular linear hypermap if the following hold:
\begin{enumerate}[{\rm(1)}]
\item $\langle r_1,r_2\rangle\cap \langle r_0,r_2\rangle=\langle r_2\rangle;$
\item $\langle r_1,r_2\rangle\langle r_0,r_2\rangle\cap \langle r_0,r_2\rangle\langle r_1,r_2\rangle=\langle r_1,r_2\rangle\cup\langle r_0,r_2\rangle.$
\end{enumerate}
The subgroups $\langle r_1,r_2\rangle, \langle r_0,r_2\rangle$
and $\langle r_0,r_1\rangle$ are stabilizers of a vertex, a hyperedge and a hyperface, respectively. Moreover, if $G=\langle r_0r_2, r_1r_2 \rangle$, then $\M$ is non-orientable; otherwise, orientable.
\end{defi}

The valency of a vertex is the number of hyperedges incident with it, and the valency of a hyperedge (a hyperface) is the number of vertices incident with it. The {\em type} of $\M$ is a triple $(k, m, n)$ of positive integers representing successively
the valencies of vertices, hyperedges and hyperfaces. If each hyperedge has valency $m$, then $\M$ is also called an $m$-{\it uniform linear hypermap}, and if each vertex has valency $k$, then $\M$ is called {\it a linear hypermap of valency $k$}. The following proposition is obvious.

\begin{prop}
The underlying linear hypergraph of a regular linear hypermap is a configuration.
\end{prop}
\begin{defi}\label{ms}
Let $\mathcal{M}=\mathcal{M}(G;r_0, r_1,r_2)$ be a linear hypermap. In this paper, we call a sequence $[g;k,m,n;V,E,F;|G|]$  a $\M$-sequence, if they stand for the
genus  of the underlying surface {\rm (also the genus of $\M$)}, the type of $\M$, the number of vertices, the number of hyperedges, the number of
hyperfaces and the number of flags  of $\M$, respectively.
\end{defi}

\begin{theo}\label{4.9}
Let $\mathcal{M}=\mathcal{M}(G;r_0, r_1,r_2)$ be a regular linear hypermap and let $H=\langle r_1,r_2 \rangle$ be the stabilizer of a vertex. Then one of the following holds.
\begin{enumerate}[{\rm(1)}]
\item $\Core_G(H)=1;$
\item $\Core_G(H)=\langle r_2\rangle$ and $G\cong G_1=\mathbb{Z}_2\times D_{2n}=\langle r_0, r_1, r_2 \di r_0^2=r_1^2=r_2^2=1, (r_0r_1)^n=(r_2r_1)^2=(r_2r_0)^2=1\rangle$ with $n\ge 3$ or $G\cong G_2 = D_{2m}= \langle r_0, r_1, r_2 \di r_0^2=r_1^2=(r_0r_1)^m=1, r_2=(r_0r_1)^{\frac{m}{2}}\rangle$ where $m$ is even and $m\ge 6$ .
\end{enumerate}
Moreover, if $2n=m$ and $n$ is odd, then the group $G_1$ is isomorphic to $G_2$. If $n$ is even, the regular linear hypermap with the automorphism group $G_1$ is isomorphic to $\mathcal{M}_1=(G_1;r_0, r_1,r_2)$; if $n$ is odd, the regular linear hypermap with the automorphism group $G_1$ is isomorphic to either $\mathcal{M}_1=(G_1;r_0, r_1,r_2)$ or $\mathcal{M}_2=(G_1;r_0, r_1r_2,r_2)$; the regular linear hypermap with the automorphism group $G_2$ is isomorphic to $\mathcal{M}_3=(G_2;r_0, r_1,r_2)$ for $4 \di m$.

If $n$ is even, then the $\mathcal{M}_1$-sequence is $[0;2,2,n;n,n,2;4n]$ and $\mathcal{M}_1$ is orientable.
If $n$ is odd, then the $\mathcal{M}_1$ and $\mathcal{M}_2$-sequences are $[0;2,2,n;n,n,2;4n]$ and $[1;2,2,2n;n,n,1;4n]$, respectively. In this case, $\mathcal{M}_2$ is non-orientable on the Projective Plane. The linear hypermap $\mathcal{M}_3$ is non-orientable on the Projective Plane with $\mathcal{M}_3$-sequence $[1;2,2,m;\frac{m}{2},\frac{m}{2},1;2m]$.
\end{theo}

\demo
Let $N=\Core_G(H)$. Then, $N\langle r_0,r_2\rangle\subseteq \langle r_1,r_2\rangle\cup\langle r_0,r_2\rangle.$ Suppose that $N\langle r_0,r_2\rangle\subsetneq\langle r_0,r_2\rangle$, then there exists $x\in N$ such that $x\langle r_0,r_2\rangle\subset\langle r_1,r_2\rangle$. As a result, $r_0\in \langle r_1,r_2\rangle$, contradicting to $\langle r_1,r_2\rangle\cap \langle r_0,r_2\rangle=\langle r_2\rangle$. So, $N\langle r_0,r_2\rangle\subseteq\langle r_0,r_2\rangle$, and then $N\le \langle r_2\rangle$.

If $N=1$, then $\Core_G(H)=1$. Note that $N$ is a normal subgroup of $G$. If $N=\langle r_2\rangle$, then either $r_2\notin \langle r_0,r_1\rangle$ or $r_2\in \langle r_0,r_1\rangle$. It follows that $$G\cong G_1=\mathbb{Z}_2\times D_{2n}=\langle r_0, r_1, r_2 \di r_0^2=r_1^2=r_2^2=1, (r_0r_1)^n=(r_2r_1)^2=(r_2r_0)^2=1\rangle,$$ where $D_{2n}$ is the dihedral group generated by $r_0, r_1$ and $n=|r_1r_0|$, or, $$G\cong G_2 = D_{2m}= \langle r_0, r_1, r_2 \di r_0^2=r_1^2=(r_0r_1)^m=1, r_2=(r_0r_1)^{\frac{m}{2}}\rangle$$ where $m$ is even. Moreover, the condition $\langle r_0,r_2\rangle \langle r_1,r_2\rangle\cap \langle r_1,r_2\rangle \langle r_0,r_2\rangle=\langle r_0,r_2\rangle\cup \langle r_1,r_2\rangle$ implies $n \ge 3$ and $m\ge 6$.
On the contrary, when $n \ge 3$ and $m\ge 6$, a simple routine check shows that $G_1$ and $G_2$ are the automorphism groups of some regular linear hypermaps.

The rest result is obvious and we omit the details.
\qed

\begin{exam}
The following Figure \ref{trH} shows two non-isomorphic regular linear hypermaps with the automorphism group $G\cong \mathbb{Z}_2\times D_{6}$. The left one is on the Sphere with $3$ hyperedges $($white components$)$
 and $2$ hyperfaces. The right one is on the Projective Plane with $3$ hyperedges $($white components$)$
 and one hyperface of valency $6$.

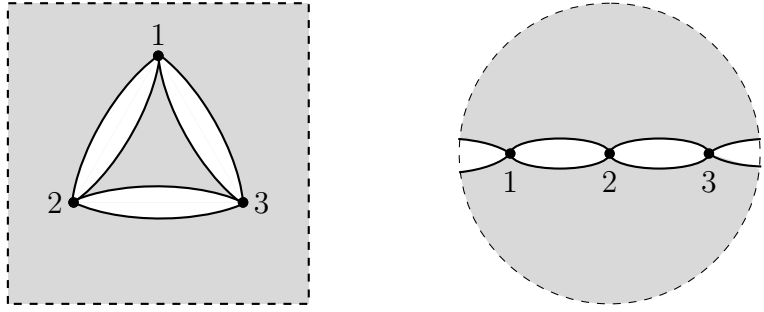
\begin{figure}[H]
\begin{center}
\begin{tikzpicture}[rotate=90]

\fill[gray!30](-2,2)--(2,2)--(2,-2)--(-2,-2)--(-2,2);
\draw[dashed, thick](-2,2)--(2,2)--(2,-2)--(-2,-2)--(-2,2);
\pgfmathsetmacro{\n}{3}
\pgfmathsetmacro{\radius}{1.3}
\foreach \i in {1,...,\n} {
\coordinate (p\i) at ({360/\n * (\i - 1)}:\radius);
\fill (p\i) circle (2pt);
}
\node at (p1) [above] {1};
\node at (p2) [left] {2};
\node at (p3) [right] {3};

\fill[white,rotate=330] (p1) arc (40:140:1.5 and 0.6);
\fill[white,rotate=150] (p2) arc (40:140:1.5 and 0.6);
\fill[white,rotate=270] (p3) arc (40:140:1.5 and 0.6);
\fill[white,rotate=90] (p2) arc (40:140:1.5 and 0.6);
\fill[white,rotate=30] (p1) arc (40:140:1.5 and 0.6);
\fill[white,rotate=210] (p3) arc (40:140:1.5 and 0.6);

\draw [thick,rotate=330] (p1) arc (40:140:1.5 and 0.6);
\draw [thick,rotate=150] (p2) arc (40:140:1.5 and 0.6);
\draw [thick,rotate=270] (p3) arc (40:140:1.5 and 0.6);
\draw [thick,rotate=90] (p2) arc (40:140:1.5 and 0.6);
\draw [thick,rotate=30] (p1) arc (40:140:1.5 and 0.6);
\draw [thick,rotate=210] (p3) arc (40:140:1.5 and 0.6);

\pgfmathsetmacro{\n}{3}
\pgfmathsetmacro{\radius}{1.3}
\foreach \i in {1,...,\n} {
\coordinate (p\i) at ({360/\n * (\i - 1)}:\radius);
\fill (p\i) circle (2pt);
}

\fill[gray!30](0,-6)circle(2);
\filldraw[fill=white, draw=black][thick] (0,-5.34) ellipse (0.2cm and 0.66cm);
\filldraw[fill=white, draw=black][thick] (0,-6.66) ellipse (0.2cm and 0.66cm);

\fill[white] (0,-4.68) arc (300:350:0.4 and 1) arc (55:115:0.45) arc (180:240:0.7 and 1);
\draw[thick] (0,-4.68) arc (300:350:0.4 and 1);
\draw[thick] (0,-4.68) arc (250:195:0.4 and 1);

\fill[white] (0,-7.32) arc (70:10:0.3 and 0.9) arc (300:240:0.45) arc (115:175:0.3 and 0.85);
\draw[thick] (0,-7.32) arc (70:10:0.3 and 0.9);
\draw[thick] (0,-7.32) arc (115:175:0.3 and 0.85);
\draw[dashed](0,-6)circle(2);

\fill (0,-6) circle (2pt);
\fill (0,-4.68) circle (2pt);
\fill (0,-7.32) circle (2pt);
\node at (-0.1,-6) [below] {2};
\node at (-0.1,-4.68) [below] {1};
\node at (-0.1,-7.32) [below] {3};

\end{tikzpicture}
\caption{Regular linear hypermaps on the Sphere (left) and Projective Plane (right) with $G\cong\mathbb{Z}_2\times D_{6}$.}
\label{trH}
\end{center}
\end{figure}

\end{exam}

We obtain the following algorithm according to Theorem \ref{4.8} and Theorem \ref{4.9}.

\begin{alm}\label{SF}
 To classify the regular linear hypermaps $\mathcal{M}$ with a given automorphism group $G$,  we need to follow two steps:
\begin{enumerate}[{\rm(1)}]
\item Check that $G$ satisfies the relations:
$$G=\langle r_0,r_1,r_2\di r_0^2=r_1^2=r_2^2=1\rangle, G_v=\langle r_1,r_2\rangle,$$
    $$\langle r_1,r_2\rangle\cap \langle r_0,r_2\rangle=\langle r_2\rangle,$$
$$\langle r_1,r_2\rangle\langle r_0,r_2\rangle\cap \langle r_0,r_2\rangle\langle r_1,r_2\rangle=\langle r_1,r_2\rangle\cup\langle r_0,r_2\rangle.$$
\item Determine the representatives of the orbits of $\Aut(G)$ on the set of triples $(r_2,r_1,r_0)$, satisfying the conditions in {\rm(1)}. Then we obtain all the non-isomorphic regular linear hypermaps $\mathcal{M}(G; r_2,r_1,r_0)$, whose number of flags is $|G|$.
\end{enumerate}
\end{alm}

 As an application of Algorithm \ref{SF}, we classify regular linear hypermaps with automorphism group isomorphic to  $A_5\times \mathbb{Z}_2$  in the following theorem.

\begin{theo}\label{A5Z2}
Suppose that $G=A_5\times \mathbb{Z}_2=A_5\times \langle a \rangle=\langle r_0, r_1, r_2\rangle, H=\langle r_1, r_2\rangle, K=\langle r_0, r_2\rangle$, where $|r_0|=|r_1|=|r_2|=|a|=2$. Then $\M(G; r_0, r_1, r_2)$ is isomorphic to one of the following linear hypermaps{\rm:}
\begin{enumerate}[{\rm(1)}]
\item $\M_1=\M(G; (12)(35),(12)(34)a,(14)(23))${\rm;}

\item $\M_2=\M(G; (12)(35),(13)(24)a,(14)(23))${\rm;}

\item $\M_3=\M(G; (12)(35)a,(12)(34)a,(13)(24)a)${\rm;}

\item $\M_4=\M(G; (12)(35)a,(13)(24)a,(12)(34)a)${\rm;}

\item $\M_5=\M(G; (12)(35)a,(14)(23),(12)(34)a)${\rm;}

\item $\M_6=\M(G; (14)(25)a,(12)(34)a,(13)(24)a)${\rm;}

\item $\M_7=\M(G; (12)(35)a,(14)(23),(13)(24)a)${\rm;}

\item $\M_8=\M(G; (14)(25),(12)(34)a,(14)(23))${\rm;}

\item $\M_9=\M(G; (14)(25)a,(14)(23),(12)(34)a)${\rm;}

\item $\M_{10}=\M(G; (14)(35)a,(12)(45),(13)(45))${\rm;}

\item $\M_{11}=\M(G; (14)(35)a,(12)(45)a,(13)(45)a)${\rm;}

\item $\M_{12}=\M(G; (13)(24)a,(12)(45)a,(13)(45)a)${\rm;}

\item $\M_{13}=\M(G; (14)(35),(12)(45)a,(13)(45)a)${\rm;}

\item $\M_{14}=\M(G; (14)(35)a,(14)(23),(13)(45))${\rm;}

\item $\M_{15}=\M(G; (15)(34)a,(14)(23),(13)(45))${\rm;}

\item $\M_{16}=\M(G; (14)(35)a,(14)(23)a,(13)(45)a)${\rm;}

\item $\M_{17}=\M(G; (15)(34)a,(14)(23)a,(13)(45)a)${\rm;}

\item $\M_{18}=\M(G; (14)(35),(14)(23)a,(13)(45)a)${\rm;}

\item $\M_{19}=\M(G; (15)(34),(14)(23)a,(13)(45)a)$.
\end{enumerate}
Moreover, all these linear hypermaps $\M_i$ for $i\in \{1,\ldots,19\}$ are non-isomorphic.
\end{theo}

\demo
By Theorem \ref{4.9}, we have $\Core_{G}(H)=1$. So $H$ is isomorphic to $D_4, D_6$ or $D_{10}$.

{\rm (1)} Suppose that $H\cong D_4$. Since $G$ has two conjugate classes of subgroups which are isomorphic to $D_4$, we may take $\langle(12)(34),(13)(24)\rangle$ and $\langle(12)(34)a,(13)(24)a\rangle$ as the representatives. By $HK\cap KH=H\cup K$, one can get
$H=\langle(12)(34)a,(13)(24)a\rangle$.
Let $T$ be a set of the representatives for the orbits of $\Aut(G)$ on the set of triples $(r_0,r_1,r_2)$, satisfying the condition $$G=\langle r_0,r_1,r_2\di r_0^2=r_1^2=r_2^2=1\rangle, G_v=H.$$
Then $$T=\{((12)(35),(12)(34)a,(14)(23)), ((12)(35),(14)(23),(12)(34)a), ((12)(35),(13)(24)a,(14)(23)), $$ $$ ((12)(35),(14)(23),(13)(24)a), ((12)(35)a,(12)(34)a,(13)(24)a), ((12)(35)a,(13)(24)a,(12)(34)a), $$ $$ ((12)(35)a,(12)(34)a,(14)(23)), ((12)(35),(14)(23),(12)(34)a), ((14)(25)a,(12)(34)a,(13)(24)a), $$ $$ ((12)(35)a,(13)(24)a,(14)(23)), ((12)(35)a,(14)(23),(13)(24)a), ((14)(25),(12)(34)a,(14)(23)), $$ $$ ((14)(25),(14)(23),(12)(34)a), ((14)(25)a,(12)(34)a,(14)(23)), ((14)(25)a,(14)(23),(12)(34)a)\}.$$
Since $K=\langle r_0, r_2\rangle$ and $HK\cap KH=H\cup K$, we have
$(r_0,r_1,r_2)\in$ $$\{((12)(35),(12)(34)a,(14)(23)), ((12)(35),(13)(24)a,(14)(23)), ((12)(35)a,(12)(34)a,(13)(24)a), $$ $$((14)(25)a,(12)(34)a,(13)(24)a), ((12)(35)a,(13)(24)a,(12)(34)a), ((12)(35)a,(14)(23),(12)(34)a), $$ $$ ((12)(35)a,(14)(23),(13)(24)a), ((14)(25),(12)(34)a,(14)(23)), ((14)(25)a,(14)(23),(12)(34)a)\}.$$
Therefore, if $H\cong D_4$, then the mutually non-isomorphic linear hypermaps are as follows:
$$\begin{array}{l}
\M_1=\M(G; (12)(35),(12)(34)a,(14)(23)){\rm;}\\
\M_2=\M(G; (12)(35),(13)(24)a,(14)(23)){\rm;}\\
\M_3=\M(G; (12)(35)a,(12)(34)a,(13)(24)a){\rm;}\\
\M_4=\M(G; (12)(35)a,(13)(24)a,(12)(34)a){\rm;}\\
\M_5=\M(G; (12)(35)a,(14)(23),(12)(34)a){\rm;}\\
\M_6=\M(G; (14)(25)a,(12)(34)a,(13)(24)a){\rm;}\\
\M_7=\M(G; (12)(35)a,(14)(23),(13)(24)a){\rm;}\\
\M_8=\M(G; (14)(25),(12)(34)a,(14)(23)){\rm;}\\
\M_9=\M(G; (14)(25)a,(14)(23),(12)(34)a).
\end{array}$$

{\rm (2)} Suppose that $H\cong D_6$. Similar to (1), the mutually non-isomorphic linear hypermaps are as follows:
$$\begin{array}{l}
\M_{10}=\M(G; (14)(35)a,(12)(45),(13)(45)){\rm;}\\
\M_{11}=\M(G; (14)(35)a,(12)(45)a,(13)(45)a){\rm;}\\
\M_{12}=\M(G; (13)(24)a,(12)(45)a,(13)(45)a){\rm;}\\
\M_{13}=\M(G; (14)(35),(12)(45)a,(13)(45)a).
\end{array}$$

{\rm (3)} Suppose that $H\cong D_{10}$. Similar to (1), the mutually non-isomorphic linear hypermaps are as follows:
$$\begin{array}{l}
\M_{14}=\M(G; (14)(35)a,(14)(23),(13)(45)){\rm;}\\
\M_{15}=\M(G; (15)(34)a,(14)(23),(13)(45)){\rm;}\\
\M_{16}=\M(G; (14)(35)a,(14)(23)a,(13)(45)a){\rm;}\\
\M_{17}=\M(G; (15)(34)a,(14)(23)a,(13)(45)a){\rm;}\\
\M_{18}=\M(G; (14)(35),(14)(23)a,(13)(45)a){\rm;}\\
\M_{19}=\M(G; (15)(34),(14)(23)a,(13)(45)a).
\end{array}$$
\qed

\begin{theo}\label{A5Z2XL}
Suppose that $\M=\M(G; r_0, r_1, r_2)$ is a regular linear hypermap, where $G=A_5\times\langle a\rangle$, $|a|=2$.
If $\M$ is orientable, then $\M$ is isomorphic to one of the following linear hypermaps{\rm :}

{\rm(1)} $\M_1=\M(G; (12)(35)a,(12)(34)a,(13)(24)a)$ with $\M$-sequence $[0;2,5,3;30,12,20;120]${\rm;}

{\rm(2)} $\M_2=\M(G; (12)(35)a,(13)(24)a,(12)(34)a)$ with $\M$-sequence $[0;2,3,5;30,20,12;120]${\rm;}

{\rm(3)} $\M_3=\M(G; (14)(25)a,(12)(34)a,(13)(24)a)$ with $\M$-sequence  $[4;2,5,5;30,12,12;120]${\rm;}

{\rm(4)} $\M_4=\M(G; (14)(35)a,(12)(45)a,(13)(45)a)$ with $\M$-sequence $[0;3,2,5;20,30,12;120]${\rm;}

{\rm(5)} $\M_5=\M(G; (13)(24)a,(12)(45)a,(13)(45)a)$ with $\M$-sequence $[5;3,3,5;20,20,12;120]${\rm;}

{\rm(6)} $\M_6=\M(G; (14)(35)a,(14)(23)a,(13)(45)a)$ with $\M$-sequence $[0;5,2,3;12,30,20;120]${\rm;}

{\rm(7)} $\M_7=\M(G; (15)(34)a,(14)(23)a,(13)(45)a)$ with $\M$-sequence $[4;5,2,5;12,30,12;120]$.

\noindent If $\M$ is non-orientable, then $\M$ is isomorphic to one of the following linear hypermaps{\rm :}

{\rm(8)} $\M_8=\M(G; (12)(35),(12)(34)a,(14)(23))$ with $\M$-sequence $[10;2,5,6;30,12,10;120]${\rm;}

{\rm(9)} $\M_9=\M(G; (12)(35),(13)(24)a,(14)(23))$ with $\M$-sequence $[14;2,5,10;30,12,6;120]${\rm;}

{\rm(10)} $\M_{10}=\M(G; (12)(35)a,(14)(23),(12)(34)a)$ with $\M$-sequence $[6;2,3,10;30,20,6;120]${\rm;}

{\rm(11)} $\M_{11}=\M(G; (12)(35)a,(14)(23),(13)(24)a)$ with $\M$-sequence $[14;2,5,10;30,12,6;120]${\rm;}

{\rm(12)} $\M_{12}=\M(G; (14)(25),(12)(34)a,(14)(23))$ with $\M$-sequence $[6;2,3,10;30,20,6;120]${\rm;}

{\rm(13)} $\M_{13}=\M(G; (14)(25)a,(14)(23),(12)(34)a)$ with $\M$-sequence $[10;2,5,6;30,12,10;120]${\rm;}

{\rm(14)} $\M_{14}=\M(G; (14)(35)a,(12)(45),(13)(45))$ with $\M$-sequence $[6;3,2,10;20,30,6;120]${\rm;}

{\rm(15)} $\M_{15}=\M(G; (14)(35),(12)(45)a,(13)(45)a)$ with $\M$-sequence $[6;3,2,10;20,30,6;120]${\rm;}

{\rm(16)} $\M_{16}=\M(G; (14)(35)a,(14)(23),(13)(45))$ with $\M$-sequence $[10;5,2,6;12,30,10;120]${\rm;}

{\rm(17)} $\M_{17}=\M(G; (15)(34)a,(14)(23),(13)(45))$ with $\M$-sequence $[14;5,2,10;12,30,6;120]${\rm;}

{\rm(18)} $\M_{18}=\M(G; (14)(35),(14)(23)a,(13)(45)a)$ with $\M$-sequence $[10;5,2,6;12,30,10;120]${\rm;}

{\rm(19)} $\M_{19}=\M(G; (15)(34),(14)(23)a,(13)(45)a)$ with $\M$-sequence $[14;5,2,10;12,30,6;120]$.
\end{theo}

\demo
{\rm (1)} Considering $\M_1=\M(G; (12)(35)a,(12)(34)a,(13)(24)a)$.Then $|r_1r_2|=|(14)(23)|\\=2$, $|r_0r_2|=|(14235)|=5$, $|r_0r_1|=|(354)|=3$ and $|G|=120$. So $V=\frac{|G|}{|\langle r_1, r_2 \rangle|}=\frac{120}{4}=30$, $E=\frac{|G|}{|\langle r_0, r_2 \rangle|}=\frac{120}{10}=12$ and $F=\frac{|G|}{|\langle r_0, r_1 \rangle|}=\frac{120}{6}=20$. Since $|\langle r_0r_2, r_1r_2 \rangle|=|\langle (14235), (14)(23)\rangle|=60$, we get $|G_1:\langle r_0r_2, r_1r_2 \rangle|=2$. So $\M_1$ is orientable. By $V+E+F-\frac{|G_1|}{2}=2-2\gamma$, we have $\gamma=0$. Hence the $\M$-sequence of $\M_1$ is $[0;2,5,3;30,12,20;120]$.

{\rm (2)--(7)} Similar with~$(1)$.

{\rm (8)} Considering $\M_8=\M(G; (12)(35),(12)(34)a,(14)(23))$. We know $|r_1r_2|=|(13)(24)a|=2$, $|r_0r_2|=|(13524)|=5$, $|r_0r_1|=|(354)a|=6$ and $|G|=120$. Then $V=\frac{|G|}{|\langle r_1, r_2 \rangle|}=\frac{120}{4}=30$, $E=\frac{|G|}{|\langle r_0, r_2 \rangle|}=\frac{120}{10}=12$ and $F=\frac{|G|}{|\langle r_0, r_1 \rangle|}=\frac{120}{12}=10$. Since $|\langle r_0r_2, r_1r_2 \rangle|=|\langle (13524), (13)(24)a\rangle|=120$, we get $|G:\langle r_0r_2, r_1r_2 \rangle|=1$. So $\M_8$ is non-orientable. From $V+E+F-\frac{|G|}{2}=2-\gamma$ we have $\gamma=10$. Hence the $\M$-sequence of $\M_8$ is $[10;2,5,6;30,12,10;120]$.

{\rm (9)--(19)} Similar with $(8)$.
\qed

\section{Platonic linear hypermaps}
According to the definition of a regular linear hypermap, the following Propositions~\ref{pro4.1}  is quite obvious, so we omit the proof.

\begin{prop}\label{pro4.1}
Let $\M=(G;r_0,r_1,r_2)$ be a regular linear hypermap. Then $\widetilde{\M}=(G;r_1,r_0,r_2)$ is a regular linear hypermap.
\end{prop}

\begin{defi}
Let $\M=(G;r_0,r_1,r_2)$ be a regular linear hypermap. Then, $\widetilde{\M}=(G;r_1,r_0,r_2)$ is called the dual linear hypermap of $\M$. If $\M$ is isomorphic to $\widetilde{\M}$, then $\M$ is called self-dual.
\end{defi}

Recall that the {\it associated map} $M$ of a linear hyermap $\M$ is a map  whose vertices are the vertices of $\M$ and whose faces include both the hyperedges and the hyperfaces of $\M$. The
map $M$ is said to be a {\it $2$-orbit} map if its automorphism group contains a subgroup that has 2 orbits on the set of all flags of $M$. A map is said to be {\it arc-transitive} if its automorphism group is arc-transitive on the underlying
graph.

\begin{prop}\label{pro4.3}
Let $\M$ be a regular linear hypermap. Then the associated map $M$ of $\M$ is  a $2$-orbit arc-transitive map.
\end{prop}
\demo
According to the definition, a hypermap automorphism is a permutation on flags preserving the incidence relations of vertices, edges and faces of the associated map. So, a hypermap automorphism induces a map automorphism. Because a flag of a hypermap contains two flags of the associated map which can not be reversed by hypermap automorphisms, there are $2$ orbits on map flags under the assumption of a regular linear hypermap.
The one-to-one correspondence of flags and arcs implies the arc transitivity of the associated graph when the linear hypermap is regular.
\qed

\begin{rem}
A $2$-orbit map may correspond to two regular linear hypermaps with different underlying linear hypergraphs.
\end{rem}

Let $M$ be a map, and let $V(M), E(M)$ and $F(M)$ be its vertex set, edge
set, and  face set, respectively. A flag of $M$ is an ordered triple $(u, e, f)$ which are mutually incident, where $u \in V, e \in E, f \in F$.  We
call $(u, f)$ a corner of $M$ if there exists an edge $e \in E$ such that $(u, e, f)$ is a flag of
$M$. Under the present assumptions concerning maps, there exist exactly two such
edges. Such pairs of edges will be said to be {\it contiguous}. The vertex set of the {\it medial map} $M^{\med}$ of a map $M$ is $E(M)$, and two vertices of
$M^{\med}$ are adjacent exactly when the corresponding edges are contiguous in
$M$, see \cite{STW}.

\begin{defi}
Let $M$ be a map whose underlying graph is simple. A linear hypermap $\M^{\med}$ is said to be  medial, if the associated map of $\M^{\med}$ is the medial map of $M$ and the set of hyperedges and hyperfaces of $\M^{\med}$ correspond in a natural way to the set $V(M)$ and $F(M)$, respectively.
\end{defi}

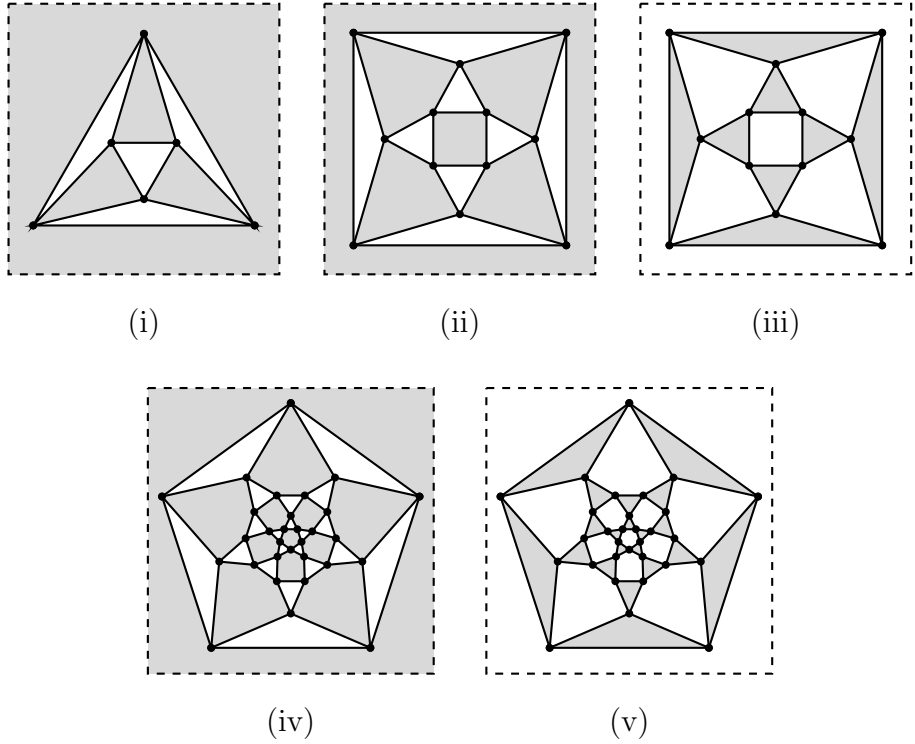
\begin{figure}[H]

\begin{center}
\begin{tikzpicture}
\begin{scope}[xshift=-4.2cm, yshift=-0.3cm]
\filldraw[fill=gray!30, draw=black, thick][dashed](-1.8,2.1)--(1.8,2.1)--(1.8,-1.5)--(-1.8,-1.5)--(-1.8,2.1);
\pgfmathsetmacro{\n}{3}
\pgfmathsetmacro{\radius}{0.5}
\foreach \i in {1,...,\n} {
\coordinate (p\i) at ({360/\n * (\i - 1)+30}:\radius);
\fill (p\i) circle (1.5pt);
}
\pgfmathsetmacro{\n}{3}
\pgfmathsetmacro{\radius}{1.7}
\foreach \i in {1,...,\n} {
\coordinate (q\i) at ({360/\n * (\i - 1)+90}:\radius);
\fill (q\i) circle (1.5pt);
}
\filldraw[fill=white, draw=black][thick](p1) -- (p2) -- (p3) -- (p1);
\filldraw[fill=white, draw=black][thick](q1) -- (p2) -- (q2) -- (q1);
\filldraw[fill=white, draw=black][thick](q1) -- (q3) -- (p1) -- (q1);
\filldraw[fill=white, draw=black][thick](p3) -- (q2) -- (q3) -- (p3);
\foreach \i in {1,...,3} {\fill (p\i) circle (1.5pt);}
\foreach \i in {1,...,3} {\fill (q\i) circle (1.5pt);}
\node at (0,-1.8) [below] {(i)};
\end{scope}
\begin{scope}
\filldraw[fill=gray!30, draw=black, thick][dashed](-1.8,1.8)--(1.8,1.8)--(1.8,-1.8)--(-1.8,-1.8)--(-1.8,1.8);
\pgfmathsetmacro{\n}{4}
\pgfmathsetmacro{\radius}{0.5}
\foreach \i in {1,...,\n} {
\coordinate (a\i) at ({360/\n * (\i - 1)+45}:\radius);
\fill (a\i) circle (1.5pt);
}
\pgfmathsetmacro{\n}{4}
\pgfmathsetmacro{\radius}{1}
\foreach \i in {1,...,\n} {
\coordinate (b\i) at ({360/\n * (\i - 1)}:\radius);
\fill (b\i) circle (1.5pt);
}
\pgfmathsetmacro{\n}{4}
\pgfmathsetmacro{\radius}{2}
\foreach \i in {1,...,\n} {
\coordinate (c\i) at ({360/\n * (\i - 1)+45}:\radius);
\fill (c\i) circle (1.5pt);
}
\filldraw[fill=white, draw=black][thick](c1) -- (c2) -- (c3) -- (c4) -- (c1);
\filldraw[fill=gray!30, draw=black][thick](a1) -- (a2) -- (a3) -- (a4) -- (a1);
\filldraw[fill=gray!30, draw=black][thick](c1) -- (b2) -- (a1) -- (b1) -- (c1);
\filldraw[fill=gray!30, draw=black][thick](c2) -- (b3) -- (a2) -- (b2) -- (c2);
\filldraw[fill=gray!30, draw=black][thick](c3) -- (b4) -- (a3) -- (b3) -- (c3);
\filldraw[fill=gray!30, draw=black][thick](c4) -- (b1) -- (a4) -- (b4) -- (c4);
\foreach \i in {1,...,4} {\fill (a\i) circle (1.5pt);}
\foreach \i in {1,...,4} {\fill (b\i) circle (1.5pt);}
\foreach \i in {1,...,4} {\fill (c\i) circle (1.5pt);}
\node at (0,-2.1) [below] {(ii)};
\end{scope}
\begin{scope}[xshift=4.2cm]
\filldraw[fill=white, draw=black, thick][dashed](-1.8,1.8)--(1.8,1.8)--(1.8,-1.8)--(-1.8,-1.8)--(-1.8,1.8);
\pgfmathsetmacro{\n}{4}
\pgfmathsetmacro{\radius}{0.5}
\foreach \i in {1,...,\n} {
\coordinate (a\i) at ({360/\n * (\i - 1)+45}:\radius);
\fill (a\i) circle (1.5pt);
}
\pgfmathsetmacro{\n}{4}
\pgfmathsetmacro{\radius}{1}
\foreach \i in {1,...,\n} {
\coordinate (b\i) at ({360/\n * (\i - 1)}:\radius);
\fill (b\i) circle (1.5pt);
}
\pgfmathsetmacro{\n}{4}
\pgfmathsetmacro{\radius}{2}
\foreach \i in {1,...,\n} {
\coordinate (c\i) at ({360/\n * (\i - 1)+45}:\radius);
\fill (c\i) circle (1.5pt);
}
\filldraw[fill=gray!30, draw=black][thick](c1) -- (c2) -- (c3) -- (c4) -- (c1);
\filldraw[fill=white, draw=black][thick](a1) -- (a2) -- (a3) -- (a4) -- (a1);
\filldraw[fill=white, draw=black][thick](c1) -- (b2) -- (a1) -- (b1) -- (c1);
\filldraw[fill=white, draw=black][thick](c2) -- (b3) -- (a2) -- (b2) -- (c2);
\filldraw[fill=white, draw=black][thick](c3) -- (b4) -- (a3) -- (b3) -- (c3);
\filldraw[fill=white, draw=black][thick](c4) -- (b1) -- (a4) -- (b4) -- (c4);
\foreach \i in {1,...,4} {\fill (a\i) circle (1.5pt);}
\foreach \i in {1,...,4} {\fill (b\i) circle (1.5pt);}
\foreach \i in {1,...,4} {\fill (c\i) circle (1.5pt);}
\node at (0,-2.1) [below] {(iii)};
\end{scope}
\end{tikzpicture}
\end{center}

\begin{center}
\begin{tikzpicture}
\begin{scope}
\filldraw[fill=gray!30, draw=black, thick][dashed](-1.9,2)--(1.9,2)--(1.9,-1.8)--(-1.9,-1.8)--(-1.9,2);
\pgfmathsetmacro{\n}{5}
\pgfmathsetmacro{\radius}{0.15}
\foreach \i in {1,...,\n} {
\coordinate (a\i) at ({360/\n * (\i - 1)-18}:\radius);
\fill (a\i) circle (1.5pt);
}
\pgfmathsetmacro{\n}{5}
\pgfmathsetmacro{\radius}{0.3}
\foreach \i in {1,...,\n} {
\coordinate (b\i) at ({360/\n * (\i - 1)+18}:\radius);
\fill (b\i) circle (1.5pt);
}
\pgfmathsetmacro{\n}{10}
\pgfmathsetmacro{\radius}{0.6}
\foreach \i in {1,...,\n} {
\coordinate (c\i) at ({360/\n * (\i - 1)}:\radius);
\fill (c\i) circle (1.5pt);
}
\pgfmathsetmacro{\n}{5}
\pgfmathsetmacro{\radius}{1}
\foreach \i in {1,...,\n} {
\coordinate (d\i) at ({360/\n * (\i - 1)-18}:\radius);
\fill (d\i) circle (1.5pt);
}
\pgfmathsetmacro{\n}{5}
\pgfmathsetmacro{\radius}{1.8}
\foreach \i in {1,...,\n} {
\coordinate (e\i) at ({360/\n * (\i - 1)+18}:\radius);
\fill (e\i) circle (1.5pt);
}
\filldraw[fill=white, draw=black][thick](e1) -- (e2) -- (e3) -- (e4) -- (e5) -- (e1);
\filldraw[fill=gray!30, draw=black][thick](e1) -- (d2) -- (c2) -- (c1) -- (d1) -- (e1);
\filldraw[fill=gray!30, draw=black][thick](e2) -- (d3) -- (c4) -- (c3) -- (d2) -- (e2);
\filldraw[fill=gray!30, draw=black][thick](e3) -- (d4) -- (c6) -- (c5) -- (d3) -- (e3);
\filldraw[fill=gray!30, draw=black][thick](e4) -- (d5) -- (c8) -- (c7) -- (d4) -- (e4);
\filldraw[fill=gray!30, draw=black][thick](e5) -- (d1) -- (c10) -- (c9) -- (d5) -- (e5);
\filldraw[fill=gray!30, draw=black][thick](a2) -- (b1) -- (c2) -- (c3) -- (b2) -- (a2);
\filldraw[fill=gray!30, draw=black][thick](a3) -- (b2) -- (c4) -- (c5) -- (b3) -- (a3);
\filldraw[fill=gray!30, draw=black][thick](a4) -- (b3) -- (c6) -- (c7) -- (b4) -- (a4);
\filldraw[fill=gray!30, draw=black][thick](a5) -- (b4) -- (c8) -- (c9) -- (b5) -- (a5);
\filldraw[fill=gray!30, draw=black][thick](a1) -- (b5) -- (c10) -- (c1) -- (b1) -- (a1);
\filldraw[fill=gray!30, draw=black][thick](a1) -- (a2) -- (a3) -- (a4) -- (a5) -- (a1);
\foreach \i in {1,...,5} {\fill (a\i) circle (1.5pt);}
\foreach \i in {1,...,5} {\fill (b\i) circle (1.5pt);}
\foreach \i in {1,...,10} {\fill (c\i) circle (1.5pt);}
\foreach \i in {1,...,5} {\fill (d\i) circle (1.5pt);}
\foreach \i in {1,...,5} {\fill (e\i) circle (1.5pt);}
\node at (0,-2.1) [below] {(iv)};
\end{scope}
\begin{scope}[xshift=4.5cm]
\filldraw[fill=white, draw=black, thick][dashed](-1.9,2)--(1.9,2)--(1.9,-1.8)--(-1.9,-1.8)--(-1.9,2);
\pgfmathsetmacro{\n}{5}
\pgfmathsetmacro{\radius}{0.15}
\foreach \i in {1,...,\n} {
\coordinate (a\i) at ({360/\n * (\i - 1)-18}:\radius);
\fill (a\i) circle (1.5pt);
}
\pgfmathsetmacro{\n}{5}
\pgfmathsetmacro{\radius}{0.3}
\foreach \i in {1,...,\n} {
\coordinate (b\i) at ({360/\n * (\i - 1)+18}:\radius);
\fill (b\i) circle (1.5pt);
}
\pgfmathsetmacro{\n}{10}
\pgfmathsetmacro{\radius}{0.6}
\foreach \i in {1,...,\n} {
\coordinate (c\i) at ({360/\n * (\i - 1)}:\radius);
\fill (c\i) circle (1.5pt);
}
\pgfmathsetmacro{\n}{5}
\pgfmathsetmacro{\radius}{1}
\foreach \i in {1,...,\n} {
\coordinate (d\i) at ({360/\n * (\i - 1)-18}:\radius);
\fill (d\i) circle (1.5pt);
}
\pgfmathsetmacro{\n}{5}
\pgfmathsetmacro{\radius}{1.8}
\foreach \i in {1,...,\n} {
\coordinate (e\i) at ({360/\n * (\i - 1)+18}:\radius);
\fill (e\i) circle (1.5pt);
}
\filldraw[fill=gray!30, draw=black][thick](e1) -- (e2) -- (e3) -- (e4) -- (e5) -- (e1);
\filldraw[fill=white, draw=black][thick](e1) -- (d2) -- (c2) -- (c1) -- (d1) -- (e1);
\filldraw[fill=white, draw=black][thick](e2) -- (d3) -- (c4) -- (c3) -- (d2) -- (e2);
\filldraw[fill=white, draw=black][thick](e3) -- (d4) -- (c6) -- (c5) -- (d3) -- (e3);
\filldraw[fill=white, draw=black][thick](e4) -- (d5) -- (c8) -- (c7) -- (d4) -- (e4);
\filldraw[fill=white, draw=black][thick](e5) -- (d1) -- (c10) -- (c9) -- (d5) -- (e5);
\filldraw[fill=white, draw=black][thick](a2) -- (b1) -- (c2) -- (c3) -- (b2) -- (a2);
\filldraw[fill=white, draw=black][thick](a3) -- (b2) -- (c4) -- (c5) -- (b3) -- (a3);
\filldraw[fill=white, draw=black][thick](a4) -- (b3) -- (c6) -- (c7) -- (b4) -- (a4);
\filldraw[fill=white, draw=black][thick](a5) -- (b4) -- (c8) -- (c9) -- (b5) -- (a5);
\filldraw[fill=white, draw=black][thick](a1) -- (b5) -- (c10) -- (c1) -- (b1) -- (a1);
\filldraw[fill=white, draw=black][thick](a1) -- (a2) -- (a3) -- (a4) -- (a5) -- (a1);
\foreach \i in {1,...,5} {\fill (a\i) circle (1.5pt);}
\foreach \i in {1,...,5} {\fill (b\i) circle (1.5pt);}
\foreach \i in {1,...,10} {\fill (c\i) circle (1.5pt);}
\foreach \i in {1,...,5} {\fill (d\i) circle (1.5pt);}
\foreach \i in {1,...,5} {\fill (e\i) circle (1.5pt);}
\node at (0,-2.1) [below] {(v)};
\end{scope}
\end{tikzpicture}
\end{center}

  \caption{Five medial linear hypermaps from Platonic solids on sphere.}\label{medial5}
\end{figure}

\begin{exam}
  In Figure {\rm \ref{medial5}}, the five medial linear hypermaps $(i)$--$(v)$ come from the five Platonic
solids -- the tetrahedron, the cube, the octahedron, the dodecahedron and
the icosahedron, respectively, described in Plato's dialogue {\it Timaeos} {\rm \cite{Si}}.  The white components and grey components depict hyperedges and hyperfaces, respectively.
\end{exam}

\begin{defi}
 Let $M$ be a map whose underlying graph is simple. The digon linear hypermap $\M^{\dig}$ of $M$ is obtained by replacing each edge of $M$ with a digon. And, the set of hyperedges and hyperfaces of $\M^{\dig}$ are the digons and the faces of $M$, respectively.
\end{defi}

\begin{theo}
  Let $M=(G;r_0,r_1,r_2)$ be a regular polyhedron whose underlying graph is simple. Then, $\M^{\dig}=(G;r_0,r_1,r_2)$ is a regular digon linear hypermap of $M$, and $\M^{\med}=(G;r_1,r_0,r_2)$ is a regular medial linear hypermap of $M$.
\end{theo}

\demo
Let $K=\lg r_2,r_0\rg$. Since $M$ is a regular polyhedron and  $H=\lg r_1, r_2\rg\cong D_{2n}$ is a dihedral group, it follows that $|K|=4$ and $K\subseteq H$ or $H\cap K=\lg r_2\rg$.

We claim that $K$ can not be contained in $H$. Otherwise, $G$ is generated by $r_1$ and $r_2$ and $n$ is even. In this case, the order of $r_0r_1$ is $2$ which implies that the width of each face of $M$ is $2$, contradicting to the assumption of a simple graph embedding. Therefore, $H\cap K=\lg r_2\rg$.

Because the underlying graph is simple, according to Lemma 7 of \cite{LS}, $H\cap r_0Hr_0=\lg r_2\rg$. As a result, $r_0H\cap Hr_0=\lg r_2\rg r_0$.
Note that $K$ is a group of order $4$. So,
$$HK=H\cup Hr_0, \hspace{.3cm} KH=r_0H\cup H$$
follows. Therefore, $HK\cap KH=K\cup H$.

Combining the preceding facts, both $ (G;r_0,r_1,r_2)$ and $ (G;r_1,r_0,r_2)$ are automorphisms of regular linear hypermaps. Considering the actions of $r_i$ and the definition of flags either in a map or in a hypermap, it is clear that $ (G;r_0,r_1,r_2)$ and $ (G;r_1,r_0,r_2)$ are automorphism groups of a regular digon linear hypermap and a regular medial linear hypermap, respectively.
\qed

\begin{cor}
    Let $\M$ be a regular linear hypermap of type $(p,q,r)$. If $p=2$ , then $\M$ is a regular medial linear hypermap; if $q=2$, then $\M$ is a regular digon linear hypermap.
\end{cor}

\begin{lem}\label{typesp}
Let $\M$ be a regular linear hypermap on the sphere. If $\M$ has $m$ flags, then, $m\geq 12$ and the $\M$-sequence is one of those listed in Table \ref{tam}.
\begin{table}[H]
\begin{center}
\caption{$\M$-sequence of regular linear hypermaps on the sphere.}
\label{tam}
\begin{tabular}{||l|l|l||}
  \hline
   $[0;3,2,3;4,6,4;24]$ & $[0;2,3,3;6,4,4;24]$ &
   $[0;3,2,4;8,12,6;48]$\\
   \hline
   $[0;2,3,4;12,8,6;48]$ &
   $[0;4,2,3;6,12,8;48]$ & $[0;2,4,3;12,6,8;48]$ \\
   \hline
   $[0;3,2,5;20,30,12;120]$ & $[0;2,3,5;30,20,12;120]$ &
   $[0;5,2,3;12,30,20;120]$\\
   \hline
   $[0;2,5,3;30,12,20;120]$ &
   $[0;2,2,\frac{m}{4};\frac{m}{4},\frac{m}{4},2;m]$ &\\
  \hline
\end{tabular}
\end{center}

\end{table}
\end{lem}

\demo
Assume that the type of $\M$ is $(p,q,r)$. Since the underlying hypergraph of $\M$ is linear and $p\ge 2, q\ge 2, r\ge 2$, according to Theorem \ref{3.2}, we get $$\frac{m}{2p}+\frac{m}{2q}+\frac{m}{2r}-\frac{m}{2}=2.$$
That is to say,
\begin{equation}\label{equ1}
1<\frac{1}{p}+\frac{1}{q}+\frac{1}{r}=\frac{4}{m}+1\leq \frac{3}{2}.
\end{equation}
So, $m\geq 8$.

If $m=8$, then $p=q=r=2$ and $\M$ has $2$ vertices, $2$ hyperedges and $2$ hyperfaces. This is contradicting to the linear assumption of $\M$. If $m=10$, then $p=q=r=5$ which contradicts to formula~\ref{equ1}. Thus, $m\geq 12$ and one can obtain Table \ref{tam}.
\qed

\begin{cor} \label{medi}
    Let $\M$ be a regular linear hypermap on the sphere. Then $\M$ is either a medial linear hypermap or a digon linear hypermap.
\end{cor}

From now on, let $G_1=A_5\times \mathbb{Z}_2=A_5\times \langle a \rangle$, $G_2=S_4$, $G_3=S_4\times \mathbb{Z}_2=S_4\times \langle a \rangle$, where $|a|=2$. By
Algorithm \ref{SF} and Theorem \ref{A5Z2XL}, we get the following theorem.
\begin{theo}
The following Table {\rm \ref{ori}} and Table {\rm \ref{non-or}} give non-isomorphic regular linear hypermaps with automorphism groups isomorphic to $G_1$, $G_2$ and $G_3$, respectively. Moreover, $\M_{1,5}$ is self-dual.
\end{theo}

\begin{table}[H]
\center
\caption{Non-isomorphic orientable regular linear hypermaps.}
\label{ori}
\begin{tabular}{||c|l|l||}
  \hline
  Groups & Regular Linear Hypermaps & $\M$-sequence\\
  \hline
  $G_1$ & $\M_{1,1}(G_1; (12)(35)a,(12)(34)a,(13)(24)a)$ & [0;2,5,3;30,12,20;120]\\
  \hline
        & $\M_{1,2}(G_1; (12)(35)a,(13)(24)a,(12)(34)a)$ & [0;2,3,5;30,20,12;120]\\
        \hline
        & $\M_{1,3}(G_1; (14)(25)a,(12)(34)a,(13)(24)a)$ & [4;2,5,5;30,12,12;120]\\
        \hline
        & $\M_{1,4}(G_1; (14)(35)a,(12)(45)a,(13)(45)a)$ & [0;3,2,5;20,30,12;120]\\
        \hline
        & $\M_{1,5}(G_1; (13)(24)a,(12)(45)a,(13)(45)a)$ & [5;3,3,5;20,20,12;120]\\
        \hline
        & $\M_{1,6}(G_1; (14)(35)a,(14)(23)a,(13)(45)a)$ & [0;5,2,3;12,30,20;120]\\
        \hline
        & $\M_{1,7}(G_1; (15)(34)a,(14)(23)a,(13)(45)a)$ & [4;5,2,5;12,30,12;120]\\
  \hline
  \hline
  $G_2$ & $\M(G_2; (13),(12),(34))$ & [0;2,3,3;6,4,4;24]\\
  \hline
        & $\M(G_2; (34),(13),(12))$ & [0;3,2,3;4,6,4;24]\\
        \hline
         \hline
  $G_3$ & $\M(G_3; (13),(12)a,(34))$ & [1;2,3,6;12,8,4;48]\\
  \hline
        & $\M(G_3; (13),(34),(12)(34)a)$ & [0;2,4,3;12,6,8;48]\\
        \hline
        & $\M(G_3; (13),(12)(34)a,(34))$ & [0;2,3,4;12,8,6;48]\\
        \hline
        & $\M(G_3; (34)a,(13),(12))$ & [1;3,2,6;8,12,4;48]\\
        \hline
        & $\M(G_3; (12)(34)a,(13),(12))$ & [0;3,2,4;8,12,6;48]\\
        \hline
        & $\M(G_3; (13)a,(12)a,(13)(24)a)$ & [0;4,2,3;6,12,8;48]\\
  \hline
\end{tabular}
\end{table}

\begin{table}[H]
\center
\caption{Non-isomorphic non-orientable regular linear hypermaps.}
\label{non-or}
\begin{tabular}{||c|l|l||}
  \hline
  Groups & Regular Linear Hypermaps & $\M$-sequence\\
  \hline
  $G_1$ & $\M(G_1; (12)(35),(12)(34)a,(14)(23))$ & [10;2,5,6;30,12,10;120]\\
  \hline
         & $\M(G_1; (14)(25)a,(14)(23),(12)(34)a)$ & [10;2,5,6;30,12,10;120]\\
        \hline
        & $\M(G_1; (12)(35),(13)(24)a,(14)(23))$ & [14;2,5,10;30,12,6;120]\\
        \hline
        & $\M(G_1; (12)(35)a,(14)(23),(13)(24)a)$ & [14;2,5,10;30,12,6;120]\\
        \hline
        & $\M(G_1; (12)(35)a,(14)(23),(12)(34)a)$ & [6;2,3,10;30,20,6;120]\\
        \hline
        & $\M(G_1; (14)(25),(12)(34)a,(14)(23))$ & [6;2,3,10;30,20,6;120]\\
        \hline
        & $\M(G_1; (14)(35)a,(12)(45),(13)(45))$ & [6;3,2,10;20,30,6;120]\\
        \hline
        & $\M(G_1; (14)(35),(12)(45)a,(13)(45)a)$ & [6;3,2,10;20,30,6;120]\\
        \hline
        & $\M(G_1; (14)(35)a,(14)(23),(13)(45))$ & [10;5,2,6;12,30,10;120]\\
        \hline
        & $\M(G_1; (14)(35),(14)(23)a,(13)(45)a)$ & [10;5,2,6;12,30,10;120]\\
        \hline
        & $\M(G_1; (15)(34)a,(14)(23),(13)(45))$ & [14;5,2,10;12,30,6;120]\\
        \hline
        & $\M(G_1; (15)(34),(14)(23)a,(13)(45)a)$ & [14;5,2,10;12,30,6;120]\\
  \hline
  \hline
  $G_2$ & $\M(G_2; (13),(12)(34),(12))$ & [1;2,3,4;6,4,3;24]\\
  \hline
        & $\M(G_2; (12)(34),(13),(12))$ & [1;3,2,4;4,6,3;24]\\
  \hline
  \hline
  $G_3$ & $\M(G_3; (13),(12)a,(12)(34)a)$ & [4;2,4,6;12,6,4;48]\\
  \hline
        & $\M(G_3; (13),(12)a,(13)(24)a)$ & [4;4,2,6;6,12,4;48]\\
  \hline
\end{tabular}
\end{table}

All vertex-transitive maps on the sphere are known. Apart from the cycles
with the two hemispheres as faces, there are 18 types and two infinite families of types (up to homeomorphisms). They correspond to the Platonic
and Archimedean solids and the families of ladders and antiprisms \cite{B,GS}.
Let $\M$ be a regular linear hypermap on the sphere. By Proposition \ref{pro4.3}, the associated map of $\M$ is arc transitive, so by Lemma \ref{typesp} and \cite{B,GS}, one can get the following theorem.
\begin{theo}\label{regsphere}
Let $ G_4=\langle r_0, r_1, r_2 \di r_0^2=r_1^2=r_0^2=1, (r_0r_1)^n=(r_1r_2)^2=(r_0r_2)^2=1\rangle \cong\mathbb{Z}_2\times D_{2n}$
with $n\ge 3$.  If $\M$ is a regular linear hypermap on the sphere, then $\M$ is isomorphic to one of the regular linear hypermaps in Table {\rm \ref{sph}}.
\end{theo}

\begin{table}[H]
\center
\caption{Non-isomorphic regular linear hypermaps on the sphere.}
\label{sph}
\begin{tabular}{||c|l|l||}
  \hline
  Groups & Regular Linear Hypermaps & $\M$-sequence\\
  \hline
  $G_1$ & $\M(G_1; (12)(35)a,(12)(34)a,(13)(24)a)$ & [0;2,5,3;30,12,20;120]\\
  \hline
        & $\M(G_1; (12)(35)a,(13)(24)a,(12)(34)a)$ & [0;2,3,5;30,20,12;120]\\
        \hline
        & $\M(G_1; (14)(35)a,(12)(45)a,(13)(45)a)$ & [0;3,2,5;20,30,12;120]\\
        \hline
        & $\M(G_1; (14)(35)a,(14)(23)a,(13)(45)a)$ & [0;5,2,3;12,30,20;120]\\
   \hline
   \hline
  $G_2$ & $\M(G_2; (13),(12),(34))$ & [0;2,3,3;6,4,4;24]\\
  \hline
        & $\M(G_2; (34),(13),(12))$ & [0;3,2,3;4,6,4;24]\\
        \hline
        \hline
  $G_3$ & $\M(G_3; (13),(34),(12)(34)a)$ & [0;2,4,3;12,6,8;48]\\
  \hline
        & $\M(G_3; (13),(12)(34)a,(34))$ & [0;2,3,4;12,8,6;48]\\
        \hline
        & $\M(G_3; (12)(34)a,(13),(12))$ & [0;3,2,4;8,12,6;48]\\
        \hline
        & $\M(G_3; (13)a,(12)a,(13)(24)a)$ & [0;4,2,3;6,12,8;48]\\
   \hline
   \hline
  $G_4$ & $\M(G_4; r_0,r_1,r_2)$ & $[0;2,2,n;n,n,2;4n]$\\
  \hline
\end{tabular}
\end{table}

%
%
%

\section{Proper regular linear hypermaps of small genus}

 Here ``proper" means that each of the parameters $p$, $q$ and $r$ in the
type $(p,q,r)$ is greater than 2.

Based on the results of \cite{C} and with the help of  {\sc Magma}\cite{BCP}, we determine the total number of proper non-isomorphic orientable regular linear hypermaps of genus 2 up to 101 in Table \ref{total}.
\begin{table}[H]
\center
\caption{Number of proper orientable regular linear hypermaps of genus 2 to 101.}
\label{total}
\begin{tabular}{||c|c|c|c|c|c|c|c|c|c|c||}
  \hline
  Genus & 2 & 3 & 4 & 5 & 6 & 7 & 8 & 9 & 10 & 11\\
  \hline
  Number & 1 & 4 & 1 & 8 & 3 & 1 & 6 & 8 & 21 & 3\\
  \hline
  \hline
  Genus & 12 & 13 & 14 & 15 & 16 & 17 & 18 & 19 & 20 & 21\\
  \hline
 Number & 0 & 14 & 2 & 11 & 5 & 19 & 0 & 21 & 0 & 9\\
  \hline
  \hline
  Genus & 22 & 23 & 24 & 25 & 26 & 27 & 28 & 29 & 30 & 31\\
  \hline
  Number & 9 & 2 & 2 & 36 & 6 & 1 & 42 & 26 & 0 & 8\\
  \hline
  \hline
  Genus & 32 & 33 & 34 & 35 & 36 & 37 & 38 & 39 & 40 & 41\\
  \hline
  Number & 0 & 62 & 0 & 0 & 5 & 67 & 0 & 0 & 13 & 35\\
  \hline
  \hline
  Genus & 42 & 43 & 44 & 45 & 46 & 47 & 48 & 49 & 50 & 51\\
  \hline
  Number & 0 & 11 & 0 & 12 & 42 & 2 & 0 & 110 & 6 & 13\\
  \hline
  \hline
  Genus & 52 & 53 & 54 & 55 & 56 & 57 & 58 & 59 & 60 & 61\\
  \hline
  Number & 0 & 0 & 0 & 61 & 8 & 51 & 1 & 2 & 0 & 36\\
  \hline
  \hline
  Genus & 62 & 63 & 64 & 65 & 66 & 67 & 68 & 69 & 70 & 71\\
  \hline
  Number & 0 & 0 & 23 & 159 & 4 & 4 & 4 & 13 & 2 & 2\\
  \hline
  \hline
  Genus & 72 & 73 & 74 & 75 & 76 & 77 & 78 & 79 & 80 & 81\\
  \hline
  Number & 0 & 156 & 0 & 0 & 43 & 2 & 5 & 4 & 0 & 88\\
  \hline
  \hline
  Genus & 82 & 83 & 84 & 85 & 86 & 87 & 88 & 89 & 90 & 91\\
  \hline
 Number & 100 & 0 & 0 & 43 & 1 & 0 & 0 & 43 & 0 & 62\\
  \hline
  \hline
  Genus & 92 & 93 & 94 & 95 & 96 & 97 & 98 & 99 & 100 & 101\\
  \hline
  Number & 6 & 5 & 7 & 0 & 4 & 279 & 0 & 18 & 33 & 17\\
  \hline
\end{tabular}
\end{table}

\section{Directions for further inquiry}
By Corollary  \ref{medi}, each regular linear hypermap on the sphere is either a medial linear hypermap or a digon linear hypermap. This observation suggests the following open problem \ref{prob}.

\begin{prob}\label{prob}
Determine the surfaces that have proper regular linear hypermaps.
\end{prob}

We classified regular linear hypermaps with automorphism group $A_5 \times \mathbb{Z}_2$ where $A_5$ is a non-abelian
simple group with smallest order in Theorem \ref{A5Z2}. This prompts the following open problem,
which is probably the most fundamental one to arise from this work.
\begin{prob}
Let $G$ be a non-abelian simple group. Classify regular linear hypermaps with automorphism group $G\times \mathbb{Z}_2$.
\end{prob}

\section*{Acknowledgement}
This work is supported by the National Natural Science Foundation of China (No. 12101535, 12231015) and the Natural Science Foundation of Beijing (No. M23017).

\end{document}